\documentclass[12pt,reqno]{amsart}

\usepackage{amsmath,amsthm,amsfonts,amscd,amssymb,euscript,graphicx,enumerate}
\numberwithin{equation}{section} 
\newcommand{\dbar}{\ensuremath{\bar \partial}}

\newcommand{\ad}{\ensuremath{\bar \partial^{*}  }}
\newcommand{\C}{\ensuremath{{\mathbb C}}}
\newcommand{\R}{\ensuremath{{\mathbb R}}}
\newcommand{\N}{\ensuremath{{\mathbb N}}}
\newcommand{\Q}{\ensuremath{{\mathbb Q}}}

\newcommand{\smooth}{\ensuremath{C^{\infty}}}
\newcommand{\smoothc}{\ensuremath{C_0 ^{\infty}}}
\newcommand{\ra}{\ensuremath{C^{\omega}}}
\newcommand{\I}{\ensuremath{\mathcal I}}
\newcommand{\II}{\ensuremath{\mathfrak I}}
\newcommand{\J}{\ensuremath{\mathcal J}}

\newcommand{\V}{\ensuremath{\mathcal V}}

\newcommand{\B}{\ensuremath{\mathfrak B}}
\newcommand{\T}{\ensuremath{\mathcal T}}
\newcommand{\ls}{\ensuremath{\mathcal L}}
\newcommand{\mt}{\ensuremath{\mathfrak M}}

\newcommand{\ct}{\ensuremath{\mathfrak C}}
\newcommand{\cv}{\ensuremath{\mathcal C}}

\newcommand{\nx}{\ensuremath{{\mathcal N}_x}}

\begin{document}
\title[Direct Proof of Termination of the Kohn Algorithm]{Direct Proof of Termination of the Kohn\\ Algorithm in the Real-Analytic Case}
\author{Andreea C. Nicoara}

\address{Department of Mathematics, University of Pennsylvania, 209 South $33^{rd}$ St.,  Philadelphia, PA 19104}

\email{anicoara@math.upenn.edu}



\begin{abstract}
In 1979 J.J. Kohn gave an indirect argument via the Diederich-Forn\ae ss Theorem showing that finite D'Angelo type implies termination of the Kohn algorithm for a pseudoconvex domain with real-analytic boundary. We give here a direct argument for this same implication using the stratification coming from Catlin's notion of a boundary system as well as algebraic geometry on the ring of real-analytic functions. We also indicate how this argument could be used in order to compute an effective lower bound for the subelliptic gain in the $\dbar$-Neumann problem in terms of the D'Angelo type, the dimension of the space, and the level of forms provided that an effective \L ojasiewicz inequality can be proven in the real-analytic case and slightly more information obtained about the behavior of the sheaves of multipliers in the Kohn algorithm.
\end{abstract}

\maketitle

\vspace{-0.37in}

\tableofcontents

\section{Introduction}

Joseph J. Kohn's solution to the $\dbar$-Neumann problem in \cite{dbneumann1} and \cite{dbneumann2}
for smooth strongly pseudoconvex domains showed subellipticity held with a gain of $\epsilon = \frac{1}{2}.$ Establishing a similar result for pseudoconvex domains proved more elusive. The breakthrough came in Kohn's 1979 Acta Mathematica paper where he had the insight of inserting a smooth function, a multiplier, in the subelliptic estimate for the $\dbar$-Neumann problem and investigating what properties the set of such multipliers had. It led to Kohn stating an algebraic algorithm for gauging the subellipticity of the $\dbar$-Neumann problem nowadays known as the Kohn algorithm. This algorithm yields sheaves of smooth germs, and since algebraic geometry on rings of $\smooth$ functions is notoriously tricky, Kohn only proved termination of this algorithm for domains with real-analytic boundary. More specifically, Kohn established the equivalence of the following three conditions for a pseudoconvex domain $\Omega$ in $\C^n$ with real-analytic boundary:

\begin{enumerate}
\item[(i)] subellipticity of the $\dbar$-Neumann problem for $(p,q)$ forms;
\item[(ii)] termination of the Kohn algorithm on $(p,q)$ forms (known as Kohn finite ideal type);
\item[(iii)] finite order of contact of holomorphic varieties of complex dimension $q$ with
the boundary of the domain $\Omega$ (finite D'Angelo $q$-type)
\end{enumerate}

After developing a fair amount of machinery, including a notion of multitype that gives more geometric information about domains satisfying condition (iii), finite D'Angelo type, David Catlin was able to prove the equivalence of conditions (i) and (iii) in a series of three papers, \cite{catlinnec}, \cite{catlinbdry}, and \cite{catlinsubell} for a smooth pseudoconvex domain
in $\C^n.$ Catlin did not, however, investigate the behavior of the Kohn algorithm as it had no obvious relationship to the machinery he had developed. For any pseudoconvex domain, implication (ii) $\rightarrow$ (i) is already a byproduct of how Kohn set up his algorithm in \cite{kohnacta}, which leaves implication (iii) $\rightarrow$ (ii) as the only one where not enough is understood. For smooth pseudoconvex domains, it is open and came to be called the Kohn Conjecture. Even for real-analytic pseudoconvex domains to which Kohn's result applies, there is no quantitative answer known, i.e. no known computation of an effective lower bound for the subelliptic gain in the $\dbar$-Neumann problem in terms of the D'Angelo type, the dimension of the space, and the level of forms. Kohn established the implication (iii) $\rightarrow$ (ii) indirectly by proving the contrapositive, namely that if the Kohn algorithm does not terminate, it means the boundary contains a real-analytic variety of holomorphic dimension at least $q,$ which by the Diederich-Forn\ae ss Theorem in \cite{df} implies a holomorphic manifold of dimension at least $q$ also sits in the boundary thus violating finite D'Angelo type.

The main result of this paper is to give a direct argument for the implication (iii) $\rightarrow$ (ii) using the stratification of the boundary into level sets of the Catlin multitype defined in \cite{catlinbdry}:

\medskip
\newtheorem{mainthm}{Main Theorem}[section]
\begin{mainthm}
\label{maintheorem} Let $\Omega$ in $\C^n$ be a pseudoconvex domain with real-analytic
boundary. Let $x_0 \in b \Omega$ be any point on the boundary of
the domain, and let the order of contact of holomorphic varieties of complex dimension $q$ with the boundary of $\Omega$ at $x_0$ be finite, i.e. $\Delta_q(b \Omega, x_0)= t < \infty.$ If $U_{x_0}$ is an appropriately small neighborhood around $x_0,$ then the Kohn algorithm on $(p,q)$ forms terminates at step $1$ densely in $U_{x_0} \cap \, b \Omega$ and by step $\min\{2n,N\}$ otherwise, where $N$ is the number of level sets of the Catlin multitype in $U_{x_0}$. $$N \leq (\lceil 2t^{n-q}\rceil -1)\lceil 2t^{n-q}\rceil^{\frac{(n-q)(n-q+1)}{2}-1},$$ where  $\lceil 2t^{n-q} \rceil$ is the ceiling of $2t^{n-q},$ i.e. the least integer greater than or equal to $2t^{n-q}.$
\end{mainthm}

The termination of the Kohn algorithm at step $1$ densely in $b \Omega$ under the assumption of finite D'Angelo type was already known in the mid 80's from work of D'Angelo in \cite{opendangelo} and Catlin in \cite{catlinbdry}. Kohn's 1979 result for real-analytic domains in \cite{kohnacta} specifies an upper bound of $2n$ for the number of steps until the Kohn algorithm terminates. In general, $N$ can be considerably larger, and yet good lower bounds for the subelliptic gain in the $\dbar$-Neumann problem follow from this construction for domains with a small number of levels of the Catlin multitype in a neighborhood compared to the dimension $n.$

The crucial idea in this work is bringing together geometric information deduced by Catlin in \cite{catlinbdry} for pseudoconvex domains of finite D'Angelo type with real algebraic geometry on the ring of real-analytic functions. Catlin's notion of boundary system yields Tougeron-Whitney elements, namely real-analytic functions with non-zero gradients that vanish on the variety corresponding to the top level of the Catlin multitype. Fortunately, the multitype is upper semi-continuous as defined, so one can inductively employ the classical \L ojasiewicz inequality for real-analytic functions proven by \L ojasiewicz in \cite{classicalloj}  to capture these Tougeron-Whitney elements inside one of the ideals of multipliers in the Kohn algorithm one level set of the Catlin multitype at a time. Once the Kohn algorithm is shown to have ended on a level set of the Catlin multitype, that information is transferred to the next step of the induction by aggregating multipliers. The latter step uses the author's result in \cite{andreeaqf} that sheaves of subelliptic multipliers corresponding to steps of the real-analytic version of the Kohn algorithm are quasi-coherent. As there are only $N$ level steps of the Catlin multitype, the Kohn algorithm will end by the $N^{th}$ step.

Some part of the argument given here is effective, so it does keep track of the lower bound for the subelliptic gain in the $\dbar$-Neumann problem in terms of the D'Angelo type, the dimension of the space, and the level of forms through the progression of the Kohn algorithm. No effective bound can be given at this time, however, due to the fact that two crucial ingredients are missing:

\begin{enumerate}
\item An effective \L ojasiewicz inequality, which we will state as a Conjectured Nullstellensatz below;
\item More information about the behavior of sheaves of subelliptic multipliers beyond the quasi-coherence proved by the author in \cite{andreeaqf}, which allows for the aggregation of subelliptic multipliers that eliminated level sets corresponding to lower values of the Catlin multitype but says nothing about the subelliptic gain corresponding to these multipliers away from the eliminated level sets.
\end{enumerate}

Given two real-analytic functions $f$ and $g$ defined on some open set $U$ of $\R^m$ and satisfying that the zero set of $g$ contains the zero set of $f,$ the classical \L ojasiewicz inequality states that for every compact subset $K$ of $U$ there exist a constant $C>0$ and some $\alpha \in \Q^+$ such that $|g(x)|^\alpha \leq C |f(x)|$ for all $x \in K;$ see \cite{classicalloj} or \cite{lojasiewicz}. For the rationality of the exponent, see \cite{bochnakrisler}. When $f$ is a real polynomial, effective Nullstellensatz results are known that compute $\alpha$ in terms of the degree of $f;$ see \cite{kollareffpol} and \cite{solerno}. When $f$ is real-analytic, no such results are known. We conjecture the following:

\medskip
\newtheorem{conjnss}[mainthm]{Conjectured Nullstellensatz}
\begin{conjnss}
Let $f$ be a $\R$-valued real-analytic function on a neighborhood $U$ of $0 \in \R^m.$ Let $x=(x_1, \dots, x_m),$ $Z=\{ x \in U \: \big| \: f(x)=0\},$ and $1 \leq \mu < m.$ Let $M \in \N^*$ be given. If the following two conditions are satisfied:
\begin{enumerate}[(a)]
\item  For every $j$ such that $1 \leq j \leq \mu,$ there exists a derivative $D^\alpha f$ such that $D^\alpha f \neq 0$ for every $x \in Z$ and the multi-index $\alpha=(\alpha_1, \dots, \alpha_m) \in \N^m$ satisfies that $\alpha_j \neq 0$ and $|\alpha| \leq M;$
\item $Z \subset \{x_j=0\}$ for every $1 \leq j \leq \mu <n.$
\end{enumerate}
Then there exist $P \in \N^*$ (computed in terms of $M$ and $m$) and positive constants $C_1, C_2, \dots, C_\mu$ such that $$|x_j|^P \leq C_j |f(x)|$$ on some potentially smaller neighborhood $V,$ $0 \in V \subset U$ for every $j$ such that $1 \leq j \leq \mu.$
\label{conjecturednullstellensatz} 
\end{conjnss}

The Conjectured Nullstellensatz~\ref{conjecturednullstellensatz} is the simplest effective \L ojasiewicz inequality that can be formulated in the real-analytic case as the role of the function $g$ in the classical statement is taken by $x_j,$ a function that defines a hypersurface and has only one non-zero derivative, the one given by $\frac{\partial}{\partial x_j}.$ This simpler Nullstellensatz suffices as it is only needed to capture Tougeron-Whitney elements coming from Catlin's boundary systems that define hypersurfaces. At first glance, condition (a) may seem over-technical, but it is actually completely natural as it postulates $f$ has a non-zero derivative transversal to $\{x_j=0\}$ of a controlled order, namely a controlled transversal vanishing order. Note that using just the vanishing order of $f$ does not work as the example $f(x,y)=x^6+y^2$ shows: $|x|^6 \leq |f(x)|,$ where the power $6$ is not related to the vanishing order of $f,$ which is $2.$ Fortunately, condition (a) can be obtained for the $(n-q)$ minors of the Levi form that kick-start the Kohn algorithm by refining the author's argument from \cite{levidet}.

In \cite{catlinsubell} Catlin obtained a lower bound $$\epsilon \geq \tau^{-n^2 \, \tau^{n^2}}$$ that holds for any smooth pseudoconvex domain in $\C^n$ and is exponential in $\tau=D_q,$ his notion of finite type. If the Conjectured Nullstellensatz~\ref{conjecturednullstellensatz} can be proven and the other sheaf theoretic issue can be sorted out, one would immediately be able to compare the effective bound yielded by our method to Catlin's effective bound via the correspondence between Catlin type and D'Angelo type derived in \cite{bazilandreea}. Other either sharp or effective bounds for subelliptic gain in the case $q=1$ were proven in \cite{siunote}, \cite{catlincho}, \cite{khanhzampieri}, and \cite{fribourgcatda} using different methods.

This paper is organized as follows: Section~\ref{kohnalg} is
devoted to the Kohn algorithm and the behavior of the sheaves it generates. Section~\ref{twotypes} surveys the results needed here that pertain to finite D'Angelo type and finite Catlin type as well as their correspondence. Section~\ref{mtbssect} introduces Catlin's boundary systems as well as his multitype and commutator multitype. Section~\ref{stratsect} then defines the notion of a Tougeron-Whitney element and links it to the Kohn algorithm via a stratification theorem. Section~\ref{effcomp} carries out effective computations of upper bounds for the number of level sets $N$ that appears in Theorem ~\ref{maintheorem} as well as for the transversal orders of vanishing of Levi minors that appear in condition (a) of the Conjectured Nullstellensatz~\ref{conjecturednullstellensatz}. Finally, Section~\ref{equivpf} proves Theorem~\ref{maintheorem} and sketches the proof of the effective version of Theorem~\ref{maintheorem} subject to the resolution of the two missing ingredients outlined above.

I am very much indebted to David Catlin, Charles Fefferman, and Pierre Milman who verified portions of this work and whose comments and suggestions greatly influenced its trajectory. I would also like to thank Francesca Acquistapace, Jason Bandlow, Vasile Brinzanescu, Fabrizio Broglia, and John D'Angelo for various useful insights.

\section{The Kohn algorithm}
\label{kohnalg}

We direct the reader to \cite{kohnacta} for full details of what we will be describing in this section. Let us begin with Kohn's definition of a subelliptic multiplier:

\smallskip
\newtheorem{subellmult}{Definition}[section]
\begin{subellmult}
Let $\Omega$ be a domain in $\C^n$ and let $x_0 \in
\overline{\Omega}.$ A $\smooth$ function $f$ is called a
subelliptic multiplier at $x_0$ for the $\bar\partial$-Neumann
problem on $\Omega$ if there exist a neighborhood $U$ of $x_0$ and
constants $C, \epsilon > 0$ such that
\begin{equation}
||\, f \varphi\, ||_\epsilon^2 \leq C \, ( \, ||\,\dbar \,
\varphi\,||^2_{\, 0} + ||\, \ad  \varphi \,||^2_{\, 0} +
||\,\varphi \,||^2_{\, 0} \,) \label{subellest}
\end{equation}
for all $(p,q)$ forms $\varphi \in \smoothc (U \cap \overline{\Omega}) \cap Dom (\ad),$ where $||\, \cdot \, ||_\epsilon$ is the Sobolev norm of
order $\epsilon$ and $||\, \cdot \,||_{\, 0}$ is the $L^2$ norm. Let $I^q (x_0)$ be the set of all subelliptic
multipliers at $x_0.$ \label{subellmultdef}
\end{subellmult}

We need two more definitions:

\medskip
\newtheorem{moduledef}[subellmult]{Definition}
\begin{moduledef}
To each $x_0 \in \overline{\Omega}$ and $q \geq 1$ we associate
the module $M^q (x_0)$ defined as the set of $(1,0)$ forms
$\sigma$ for which there exist a neighborhood $U$ of $x_0$
and constants $C, \epsilon > 0$ such that
\begin{equation}
||\, int(\bar\sigma)\, \varphi\, ||_{\, \epsilon}^2 \leq C \, ( \,
||\,\dbar \, \varphi\,||^2_{\, 0} + ||\, \ad  \varphi \,||^2_{\,
0} + ||\,\varphi \,||^2_{\, 0} \,) \label{moduleest}
\end{equation}
for all $(p,q)$ forms $\varphi \in \smoothc (U) \cap Dom (\ad),$
where $int(\bar\sigma)\varphi$ denotes the interior multiplication
of the $(0,1)$ form $\bar\sigma$ with the $(p,q)$ form $\varphi.$
\end{moduledef}

\medskip
\newtheorem{realraddef}[subellmult]{Definition}
\begin{realraddef}
Let $J \subset \smooth (x_0),$ the ring of germs of smooth functions at $x_0,$ then the real radical of $J$
denoted by $\sqrt[\R]{J}$ is the set of $g \in \smooth (x_0)$ such
that there exists some $f \in J$ and some positive natural number
$m \in \N^\ast$ such that $$|g|^m \leq |f|$$ on some neighborhood
of $x_0.$
\end{realraddef}

\medskip The $\dbar$-Neumann problem is subelliptic inside the domain $\Omega,$ so we are only interested in the case $x_0 \in b \Omega.$ Theorem $1.21$ of \cite{kohnacta} gives the properties of $I^q(x_0):$

\medskip
\newtheorem{subellcor}[subellmult]{Corollary}
\begin{subellcor}
If $\Omega$ is a smooth pseudoconvex domain and if $x_0
\in \overline{\Omega},$ then we have:
\begin{enumerate}
\item[(a)] $I^q (x_0)$ is an ideal.
\item[(b)] $I^q (x_0) = \sqrt[\R]{I^q (x_0)}.$
\item[(c)] If $r=0$ on $b \Omega,$ then $r \in I^q (x_0)$ and the
coefficients of $\partial r \wedge \dbar r \wedge (\partial \dbar
r)^{n-q}$ are in $I^q (x_0).$
\item[(d)] If $f_1, \dots, f_j \in I^q (x_0),$ then the
coefficients of $ \partial r \wedge \dbar r \wedge \partial f_1 \wedge \dots \wedge \partial f_j
\wedge (\partial \dbar
r)^{n-q-j}$ are in $I^q (x_0)$ for $j \leq
n-q.$\label{subellpropcor}
\end{enumerate}
\end{subellcor}

\bigskip\noindent Examining the proof of Theorem 4.7 in Section 4 of \cite{kohnacta} yields all the necessary information about the cost in terms of the reduction in the subelliptic gain in the $\dbar$-Neumann problem of performing each operation that gives rise to new subelliptic multipliers:

\medskip
\newtheorem{subellcost}[subellmult]{Proposition}
\begin{subellcost}
If $\Omega$ is a smooth pseudoconvex domain and if $x_0 \in \overline{\Omega},$
then $I^q (x_0)$ and $M^q (x_0)$ have the following properties:
\begin{enumerate}
\item[(i)] If $x_0 \in  \overline{\Omega},$ then $r \in I^q (x_0)$ satisfies \eqref{subellest} with $\epsilon=1.$
\item[(ii)] If $x_0 \in b \Omega$ and $\theta$ is any smooth $(0,1)$ form such
that $\langle \theta, \dbar r \rangle = 0$ on $b \Omega,$ then $int (\theta)\, \partial
\dbar r \in M^q (x_0)$ satisfies \eqref{moduleest} with $\epsilon=\frac{1}{2}.$
\item[(iii)] If $f \in I^q (x_0)$ satisfies \eqref{subellest} with some $\epsilon>0$ and if $g \in \smooth (x_0)$ is such that $|g| \leq |f|$ in a neighborhood of $x_0,$ then $g \in I^q (x_0)$ satisfies \eqref{subellest} with the exact same $\epsilon.$
\item[(iv)]  If $f \in I^q (x_0)$ satisfies \eqref{subellest} with some $\epsilon>0$ and if $g \in \smooth (x_0)$ is such that $|g|^m \leq |f|$ for an integer $m \in \N^\ast$ in a neighborhood of $x_0,$ then $g \in I^q (x_0)$ satisfies \eqref{subellest} with $\frac{\epsilon}{m}.$
\item[(v)] If $f \in I^q (x_0)$ satisfies \eqref{subellest} with some $\epsilon>0,$ then $\partial f \in M^q (x_0)$ satisfies \eqref{moduleest} with $\frac{\epsilon}{2},$ where $\partial f$ is the complex gradient of $f.$
\item[(vi)] If $\sigma_1, \dots, \sigma_{n+1-q} \in M^q (x_0)$ satisfy \eqref{moduleest} with $\epsilon_1, \dots, \epsilon_{n+1-q}$ respectively, then the coefficients of their wedge product $\det_{n-q+1} (\sigma_1, \dots, \sigma_{n+1-q}) \in I^q (x_0)$ satisfy \eqref{subellest} with $\epsilon=\min_{1 \leq j \leq n+1-q}\, \epsilon_j.$ \label{subellcostprop}
\end{enumerate}
\end{subellcost}

\smallskip\noindent {\bf Remarks:} 

\smallskip\noindent (1) Pseudoconvexity is essential for (ii) and (vi) and irrelevant for the rest.

\smallskip\noindent (2) None of the operations in this proposition require a shrinking of neighborhood as it can be seen by carefully reading Section 4 of \cite{kohnacta}, but the neighborhood needs to be small enough that special frames of vector fields and dual forms can be defined. We will describe this construction below.

\medskip\noindent {\bf The Kohn Algorithm:}

\medskip\noindent {\bf Step 1}  $$I^q_1(x_0) = \sqrt[\R]{(\, r,\,
\text{coeff}\{\partial r \wedge \dbar r \wedge (\partial \dbar
r)^{n-q}\}\, )}$$

\medskip\noindent {\bf Step (k+1)} $$I^q_{k+1} (x_0) = \sqrt[\R]{(\, I^q_k
(x_0),\, A^q_k (x_0)\, )},$$ where $$A^q_k (x_0)= \text{coeff}\{\partial
f_1 \wedge \dots \wedge \partial f_j \wedge \partial r \wedge
\dbar r \wedge (\partial \dbar r)^{n-q-j}\}$$ for $f_1, \dots, f_j
\in I^q_k (x_0)$ subelliptic multipliers at step $k$ and $j \leq n-q.$ As usual $( \, \cdot \, )$
stands for the ideal generated in the ring $\smooth (x_0)$ and $coeff\{\partial r \wedge \dbar r \wedge (\partial\dbar r)^{n-q}\}$ is the determinant of the Levi form for $q=1,$ whereas for $q>1$ it consists of all $(n-q)$ minors of the Levi form. The algorithm terminates when a unit is captured inside $I^q_k(x_0).$

In the setting of Theorem~\ref{maintheorem}, the domain $\Omega$ we consider is defined by a real-analytic function $r.$ Therefore, just like Kohn does in section 6 of \cite{kohnacta}, we can consider the following modification of the algorithm: $$\tilde I^q_1(x_0) = \sqrt[\R]{(\, r,\,coeff\{\partial r \wedge \dbar r \wedge (\partial \dbar r)^{n-q}\}\, )_{\ra(x_0)}}$$ and $$\tilde I^q_{k+1} (x_0) = \sqrt[\R]{(\, \tilde I^q_k (x_0),\, \tilde A^q_k (x_0)\, )_{\ra(x_0)}},$$ where $$\tilde A^q_k (x_0)= coeff\{\partial
f_1 \wedge \dots \wedge \partial f_j \wedge \partial r \wedge
\dbar r \wedge (\partial \dbar r)^{n-q-j}\}$$ for $f_1, \dots, f_j
\in\tilde I^q_k (x_0)$ and $j \leq n-q.$ Here $\ra(x_0)$ is the ring of real-analytic germs at $x_0 \in b\Omega,$ and the algebraic operations of the algorithm, generating ideals and taking real radicals, are taking place just in $\ra(x_0)$ as the subscript indicates. Obviously, the termination of the modified algorithm implies the termination of the original algorithm. By its very definition, the modified Kohn algorithm generates an increasing chain of ideals
$$\tilde I^q_1 (x_0) \subset \tilde I^q_2 (x_0) \subset \cdots$$ in the Noetherian ring $\ra (x_0),$ so we know this chain of ideals stabilizes. Our task is showing it stabilizes at the ring itself. At certain points of the argument, we will consider the behavior of the algorithm on a neighborhood $U \ni x_0,$ i.e. we will be examining $\tilde I^q_k (U).$ To relate the behavior of $\tilde I^q_k(x_0)$ with that of $\tilde I^q_k(U),$ we need to understand the properties of sheaves of subelliptic multipliers. For all $k \geq 1,$ we denote by $\tilde \II^q_k$ the sheaf of real-analytic subelliptic multipliers obtained at step $k$ of this modification of the Kohn algorithm and by $\tilde \II^q$ the sheaf of real-analytic subelliptic multipliers for the $\dbar$-Neumann problem on $(p,q)$ forms. We recall from \cite{andreeaqf} the main result governing the behavior of the sheaves $\tilde \II^q_k$ proved via the concept of a quasi-flasque sheaf that Jean-Claude Tougeron defined in \cite{tougeronqf}: 

\medskip
\newtheorem{qfsteps}[subellmult]{Theorem}
\begin{qfsteps}
Let $\Omega$ in $\C^n$ be a domain with real-analytic boundary $b \Omega.$ Let $\tilde U$ be any open subset of $b \Omega$ such that $\tilde U$ is contained in a compact semianalytic subset $Y$ of $b \Omega.$ If $\Omega$ is bounded, $b \Omega$ itself may be taken as $\tilde U.$ The ideal sheaf $\tilde \II^q$ of real-analytic subelliptic multipliers for the $\dbar$-Neumann problem on $(p,q)$ forms defined on $\tilde U$ is coherent. Additionally, if $\Omega$ is pseudoconvex, the multiplier ideal sheaf $\tilde \II^q_k$ given by the modified Kohn algorithm on $\tilde U$ at step $k$ for each $k \geq 1$ is also coherent. In other words, $\tilde \II^q$ and $\tilde \II^q_k$ for all $k \geq 1$ are quasi-coherent sheaves.
\label{quasiflasquesteps}
\end{qfsteps}

As promised in the second remark after Proposition~\ref{subellcostprop}, we now recall the standard type of neighborhood used by Kohn in \cite{kohnacta} except that we exchange indices $1$ and $n$ in order to be consistent with \cite{catlinbdry} whose boundary system construction comes into our argument later on. We choose a defining function $r$ for the domain $\Omega$ such
that $|\partial r|_x = 1$ for all $x$ in a neighborhood of $b
\Omega.$ Let $U$ be a neighborhood of $x_0$ small enough that
the previous condition holds on $U.$ We choose $(1,0)$ forms
$\omega_1, \dots, \omega_n$ on $U$ satisfying that $\omega_1 =
\partial r$ and $\langle \omega_i, \omega_j \rangle = \delta_{ij}$
for all $x \in U.$ We define by duality $(1,0)$ vector fields
$L_1, \dots, L_n$ such that $\langle \omega_i, L_j \rangle =
\delta_{ij}$ for all $x \in U.$ Thus, on $U \cap b
\Omega,$ $$L_j (r) = \bar L_j (r) = \delta_{1j}.$$ We define a
vector field $T$ on $U \cap b \Omega$ by $$T = L_1 - \bar L_1.$$
The collection of vector fields $L_2, \dots, L_n,
\bar L_2, \dots, \bar L_n, T$ gives a local basis for the
tangent space $T (U \cap b \Omega).$ A $(p,q)$
form $\varphi$ can be expressed in terms of the corresponding
local basis of dual forms on $U$ as $$\varphi = \sum_{|I|=p, \:
|J|=q} \: \varphi_{IJ} \, d\omega_I \wedge d \bar \omega_J,$$ for $I$
and $J$ multi-indices in $\N^n.$ As Kohn shows in \cite{kohnacta}, $\varphi \in Dom
(\ad)$ means precisely that
$$\varphi_{IJ} (x) = 0$$ when $1 \in J$ and $x \in b \Omega.$ The
Levi form is also computed in this local basis. \label{prenbd}

The neighborhood $U$ described above is not yet the same as the
neighborhood $U_{x_0}$ that appears in the statement of the Main
Theorem~\ref{maintheorem}. Two other conditions we
will impose later on that may shrink $U$ further, one condition will ensure Theorem~\ref{quasiflasquesteps} holds on $U_{x_0}$ and the other condition will force D'Angelo type to be finite and effectively bounded on all of $U_{x_0}.$ The latter will be discussed in Section~\ref{twotypes}.

\medskip Following \cite{kohnacta} let us now define the notion of Zariski tangent space to an
ideal and to a variety, which will allow us to introduce the notion of holomorphic dimension of a variety. We will then recall from \cite{catlinbdry} Catlin's definition of the holomorphic dimension of a variety, which is slightly different from Kohn's. These concepts will be used in Section~\ref{stratsect}.

\smallskip
\newtheorem{zariskidef}[subellmult]{Definition}
\begin{zariskidef}
Let $\I$ be an ideal in $\smooth (U)$ and let $\V (\I)$ be the
variety corresponding to $\I.$ If $x \in \V (\I),$ then we define
$Z^{\, 1,0}_x (\I)$ the Zariski tangent space of $\I$ at $x$ to be
$$Z^{\, 1,0}_x (\I) = \{ \, L \in T^{\, 1,0}_x (U) \: | \: L(f)=0 \:\:
\forall \: f \in \I \, \},$$ where $T^{\, 1,0}_x (U)$ is the
$(1,0)$ tangent space to $U \subset \C^n$ at $x.$ If $\V$ is a
variety, then
$$Z^{\, 1,0}_x (\V) = Z^{\, 1,0}_x (\I (\V)), $$ where $\I (\V)$
is the ideal of all functions in $\smooth (U)$ vanishing on $\V.$
\label{zariskidefinition}
\end{zariskidef}

\medskip\noindent The next lemma is Lemma $6.10$ of \cite{kohnacta}
that relates $Z^{\, 1,0}_x (\I)$ with $Z^{\, 1,0}_x (\V(\I)):$

\smallskip
\newtheorem{zariskilemma}[subellmult]{Lemma}
\begin{zariskilemma}
\label{zariski} If $\I$ is an ideal in $\smooth (U)$ and $x \in
\V (\I),$ then
\begin{equation}
Z^{\, 1,0}_x (\V(\I)) \subset Z^{\,1,0}_x (\I). \label{zarexpr}
\end{equation}
Equality holds in ~\eqref{zarexpr} if the ideal $\I$ satisfies the
Nullstellensatz, namely $\I = \I (\V(\I)).$
\end{zariskilemma}

\medskip\noindent If $\I$ is an ideal in $\ra(U),$ then $\I (\V(\I))$ is computed in $\ra(U)$ for the purposes of both Definition~\ref{zariskidefinition} and Lemma~\ref{zariski}. Let $$\nx = \{ \, L \in T^{\, 1,0}_x (b \Omega) \: |
\: \langle \, (\partial \dbar r)_x \, , \, L \wedge \bar L \,
\rangle = 0 \, \}.$$ $\nx$ is the subspace of $T^{\, 1,0}_x (b \Omega)$
consisting of the directions in which the Levi form vanishes. We
can now give the definition of the holomorphic dimension of a
variety sitting in the boundary of the domain $\Omega$ first according to Kohn in \cite{kohnacta} and then according to Catlin in \cite{catlinbdry}:

\smallskip
\newtheorem{holdimdefk}[subellmult]{Definition (Kohn)}
\begin{holdimdefk}
Let $\V$ be a variety in $U$ that corresponds to an ideal $\I$ in
$\smooth (U)$ or $\ra(U)$ and satisfies $\V \subset b \Omega.$ We define the
holomorphic dimension of $\V$ in the sense of Kohn by $$hol. \: dim \, (\V) = \min_{x
\in \V} \dim Z^{\, 1,0}_x (\V) \cap \nx.$$ \label{holdimkohn}
\end{holdimdefk}

\smallskip
\newtheorem{holdimdefc}[subellmult]{Definition (Catlin)}
\begin{holdimdefc}
Let $\V$ be a variety in $U$ that corresponds to an ideal $\I$ in
$\smooth (U)$ or $\ra(U)$ and satisfies $\V \subset b \Omega.$ We define the
holomorphic dimension of $\V$ in the sense of Catlin by $$hol. \: dim \, (\V) = \max_{x
\in \V} \dim Z^{\, 1,0}_x (\V) \cap \nx.$$ \label{holdimcatlin}
\end{holdimdefc}

\medskip\noindent We can now state the Diederich-Forn\ae ss Theorem in \cite{df} mentioned in the introduction:

\smallskip
\newtheorem{dfthm}[subellmult]{Theorem}
\begin{dfthm}
Let $W$ be a smooth pseudoconvex real-analytic hypersurface in $\C^n.$ Suppose $S \subset W$ is a not necessarily closed real-analytic subvariety with $hol. \: dim \, (S) =q.$ Let $z_0 \in S$ be an arbitrary point and $U = U(z_0)$ an open neighborhood of $z_0.$ Then there exists a complex submanifold $V \subset U \cap W$ of dimension at least $q.$ The manifold $V$ can always be chosen in such a way that $S \cap V \neq \emptyset$ and that in fact $hol. \: dim \, (S \cap V) = q.$ \label{dftheorem} 
\end{dfthm}

\smallskip\noindent {\bf Remark:} In the Diederich-Forn\ae ss Theorem, the holomorphic dimension is meant in the sense of Kohn.

\medskip\noindent Kohn's Proposition $6.12$ of \cite{kohnacta} provides an equivalent condition to the Kohn algorithm not advancing at a particular point:

\smallskip
\newtheorem{stillthere}[subellmult]{Proposition}
\begin{stillthere}
If $x \in \V(I^q_k (x_0)),$ then $$x \in \V(I^q_{k+1} (x_0))  \quad
\Leftrightarrow \quad \dim (Z^{\, 1,0}_x (I^q_k (x_0)) \cap \nx)
\geq q.$$ \label{stillthereprop}
\end{stillthere}

\smallskip\noindent {\bf Remark:} Kohn's proof of this proposition also applies to the real-analytic version of the Kohn algorithm, i.e. if $x \in \V(\tilde I^q_k (x_0)),$ then $$x \in \V(\tilde I^q_{k+1} (x_0))  \quad
\Leftrightarrow \quad \dim (Z^{\, 1,0}_x (\tilde I^q_k (x_0)) \cap \nx)
\geq q.$$

\medskip\noindent Kohn uses this result in a fundamental way in \cite{kohnacta} in order to show that if the Kohn algorithm does not terminate, then the boundary $b \Omega$ cannot have finite D'Angelo type at $x_0.$ By contrast, we will simply point out in the proof of Theorem~\ref{maintheorem} where this proposition could be used and show that our method yields more information. We now close the section with the \L ojasiewicz Nullstellensatz from \cite{lojasiewicz}:

\smallskip
\newtheorem{lojnss}[subellmult]{Theorem }
\begin{lojnss}
If $\J$ is an ideal of $\ra(x_0),$ then $\I (\V(\J))=\sqrt[\R]{\J}.$
\label{lojasiewicznss}
\end{lojnss}

\smallskip\noindent {\bf Remark:} This result obviously follows from the classical \L ojasiewicz inequality stated in the introduction.

\section{Notions of finite type}
\label{twotypes}

We will define finite D'Angelo type $\Delta_q$ here but not finite Catlin type $D_q.$ The reader is directed to \cite{opendangelo} and \cite{dangelo} for comprehensive discussions of $\Delta_q,$ to \cite{catlinsubell} where $D_q$ is introduced, and to \cite{bazilandreea} that relates the two.

Let $\cv=\cv(n,x_0)$ be the set of all germs of holomorphic curves $$\varphi: (U,0) \rightarrow (\C^n, x_0),$$ where $U$ is some neighborhood of the origin in $\C^1$ and $\varphi(0)=x_0.$ Let ${\text ord}_0 \, \varphi_j$ be the order of vanishing of the $j^{th}$ component of $\varphi$ at $0.$ We set ${\text ord}_0 \, \varphi = \min_{1 \leq j \leq n} \, {\text ord}_0 \, \varphi_j.$ 

\medskip
\newtheorem{firstft}{Definition}[section]
\begin{firstft}
Let $W$ be a real hypersurface in $\C^n$ and $r$ a defining function for $W.$ The D'Angelo $1$-type at $x_0 \in W$ is given by $$\Delta_1 (W, x_0) = \sup_{\varphi \in \cv(n,x_0)} \frac{ {\text ord}_0 \, \varphi^* r}{{\text ord}_0 \, \varphi},$$ where $\varphi^* r$ is the pullback of $r$ under $\varphi.$ If $\Delta_1 (W, x_0)$ is finite, we call $x_0$ a point of finite D'Angelo $1$-type.
\end{firstft}

\medskip
\newtheorem{qft}[firstft]{Definition}
\begin{qft}
\label{qfinitetype}
Let $W$ be a real hypersurface in $\C^n$ and $r$ a defining function for $W.$ The D'Angelo $q$-type at $x_0 \in W$ for $q>1$ is given by $$\Delta_q (W, x_0) =\inf_\phi \sup_{\varphi \in \cv(n-q+1,x_0)} \frac{ {\text ord}_0 \, \varphi^* \phi^*r}{{\text ord}_0 \, \varphi}=\inf_\phi \Delta_1 (\phi^*r, x_0),$$ where $\phi : \C^{n-q+1} \rightarrow \C^n$ is any linear embedding of $\C^{n-q+1}$ into $\C^n$ and we have identified $x_0$ with $\phi^{-1} (x_0).$ If $\Delta_q (W, x_0)$ is finite, we call $x_0$ a point of finite D'Angelo $q$-type.
\end{qft}

We will concentrate now just on the results involving $\Delta_q$ and $D_q$ that are essential for our argument here. We start with Theorem 6.2 from p.634 of \cite{opendangelo}:

\medskip
\newtheorem{deltaqopen}[firstft]{Theorem}
\begin{deltaqopen}
Let $W$ be a smooth real hypersurface in $\C^n.$
Let $\Delta_q (W, x_0)$ be finite at some $x_0 \in W,$ then there exists a neighborhood $V$ of $x_0$ on which $$\Delta_q (W, x) \leq 2 (\Delta_q(W, x_0))^{n-q}.$$ \label{deltaqopenthm}
\end{deltaqopen}

\medskip\noindent The next result from \cite{bazilandreea} puts $\Delta_q$ in correspondence to $D_q:$

\newtheorem{correspt}[firstft]{Theorem}
\begin{correspt}
\label{typescorrespondence}
 Let $\Omega$ in $\C^n$ be a domain with $\smooth$
boundary. Let $x_0 \in b \Omega$ be a point on the boundary of
the domain, and let $1 \leq q <n.$
 If $b \Omega$ is pseudoconvex at $x_0$ and $\Delta_q(b \Omega, x_0)<\infty,$ then $$D_q(b \Omega, x_0) \leq \Delta_q(b \Omega, x_0) \leq 2 \left(\frac{ D_q(b \Omega, x_0)}{2} \right)^{n-q}.$$
\end{correspt}

\smallskip
\newtheorem{nbd}[firstft]{Properties of the neighborhood $U_{x_0}$}
\begin{nbd}
The neighborhood $U_{x_0} \ni x_0$ in which we will run the Kohn
algorithm satisfies: \label{nbdprop}
\begin{enumerate}
\item $U_{x_0} \subset U,$ where $U$ is the neighborhood of $x_0$ described on
page~\pageref{prenbd} of Section~\ref{kohnalg};
\item Theorem~\ref{quasiflasquesteps} holds on $U_{x_0};$
\item $b\Omega \cap U_{x_0}$ satisfies Theorem~\ref{deltaqopenthm};
\item The closure $\overline{U}_{x_0}$ is compact in $\C^n.$
\end{enumerate}
\end{nbd}

\section{Catlin's multitype and boundary systems}
\label{mtbssect}

We will briefly recall here Catlin's concepts of boundary system, multitype $\mt(x_0),$ and commutator multitype $\ct(x_0)$ from \cite{catlinbdry}. The reader is directed to \cite{catlinbdry} or \cite{levidet} for more details. $\mt(x_0)$ and $\ct(x_0)$ are $n$-tuples of rational numbers satisfying certain properties. Catlin called all such $n$-tuples weights.

\medskip
\newtheorem{wgt}{Definition}[section]
\begin{wgt}
\label{wgtdef} Let $\Gamma_n$ denote the set of $n$-tuples of rational numbers $\Lambda=(\lambda_1, \dots, \lambda_n)$ with $1 \leq \lambda_i \leq + \infty$ satisfying the following two properties:
\begin{enumerate}
\item[(i)] $\lambda_1 \leq \lambda_2 \leq \cdots \leq \lambda_n.$
\item[(ii)] For each $k$ such that $1 \leq k \leq n,$ either $\lambda_k = + \infty$ or there exists a set of integers $a_1, \dots , a_k$ such that $a_j >0$ for all $1 \leq j \leq k$ and $$\sum_{j=1}^k \frac{a_j}{\lambda_j}=1.$$
\end{enumerate}
The set $\Gamma_n$ is ordered lexicographically, i.e. given $\Lambda', \Lambda'' \in \Gamma_n$ such that $\Lambda'=(\lambda'_1, \dots, \lambda'_n)$ and $\Lambda''=(\lambda''_1, \dots, \lambda''_n),$ then $\Lambda' < \Lambda''$ if there exists $k$ with $1 \leq k \leq n$ such that $\lambda'_j = \lambda''_j$ for all $j<k$ and $\lambda'_k < \lambda''_k.$ The set $\Gamma_n$ is called the set of weights.
\end{wgt}

\smallskip\noindent {\bf Remark:} At times we may work with truncated weights, which are $\nu$-tuples of rational numbers $(\lambda_1, \dots, \lambda_\nu)$ satisfying this definition. We denote by $\Gamma_\nu$ the set of such weights.

\medskip\noindent Let $\Gamma'$ be a set of weights. A weight $\Lambda$ is said to dominate all weights in $\Gamma'$ if $\Lambda \geq \Lambda'$ for every $\Lambda' \in \Gamma'.$ Catlin called distinguished weights all elements of $\Gamma_n$ that dominate the vanishing of the defining function $r$ of the domain $\Omega$ in various directions. The multitype $\mt(x_0)$ is then the smallest weight in $\Gamma_n$ that dominates all the distinguished weights. Remarkably enough, Catlin was able to construct another weight $\ct(x_0),$ the commutator multitype, by differentiating the defining function $r$ in a controlled manner and to show $\mt(x_0)=\ct(x_0)$ when the domain $\Omega$ is pseudoconvex. The notion of boundary system $\B_\nu (x_0)$ is a byproduct of the construction of $\ct(x_0).$

The commutator multitype  $\ct(x_0)=(c_1, \dots, c_n) \in \Gamma_n$ always satisfies that $c_1=1$ because as explained on page~\pageref{prenbd}, $L_1(r)=1.$ Set $r_1=r.$ Let $p$ be the rank of the Levi form of $b \Omega$ at $x_0.$ Set $c_i=2$ for $i=2, \dots, p+1.$  In the construction on page~\pageref{prenbd}, choose the smooth vector fields of type $(1,0)$ $L_2, \dots, L_{p+1}$ such that $L_i(r)=\partial r (L_i)\equiv 0$ and the $p \times p$ Hermitian matrix $\partial \bar \partial r (L_i, L_j)(x_0)$ for $2 \leq i,j \leq p+1$ is nonsingular. Round parentheses denote the evaluation of forms on vector fields. If $p+1 \geq \nu,$ $\ct^\nu (x_0)=(1,2, \dots,2),$ and we are done.

If $p+1 < \nu,$ consider the $(1,0)$ smooth vector fields in the kernel of the Levi form at $x_0.$ Let $T^{\, (1,0)}_{p+2}$ be the bundle consisting of $(1,0)$ vector fields $L$ such that $\partial r (L)=0$ and $\partial \bar \partial r (L, \bar L_j)=0$ for $j=2, \dots, p+1,$ and let $\T_{p+2}$ be the set of germs of sections of $T^{\, (1,0)}_{p+2}.$ It is obvious we now have to consider lists of vector fields of length at least $3$ in order to pinpoint other directions in which the defining function vanishes to finite order besides those involved in the non-singular part of the Levi form. Let $l \in \N$ be such that $l \geq 3,$ and let $\ls$ be a list of vector fields $\ls= \{L^1, \dots, L^l\}$ with $L^j=L$ or $\bar L^j = \bar L$ for every $j,$ $1 \leq j \leq l,$ where $L \in T^{\, (1,0)}_{p+2}$ is a fixed, non-vanishing vector field. Let $\ls  \partial r$ be the function $$\ls  \partial r (x) = L^1 \cdots L^{l-2} \, \partial r \, ([L^{l-1},L^l]) (x)$$ for $x \in b \Omega.$ Note that one of $L^{l-1}$ and $L^l$ needs to be a $(1,0)$ vector field and the other one a $(0,1)$ vector field for their bracket not to be identically zero. If $\ls \partial r (x_0)=0$ for every such list $\ls,$ we set $c_{p+2}= \infty;$ otherwise, there exists at least one list $\ls$ such that $\ls \partial r (x_0)\neq 0.$ In the latter case, choose $\ls$ of minimal length $l$ and set $c_{p+2}=l.$ Note that $L=X+iY$ for $X$ and $Y$ $\R$-valued vector fields, and define functions $$f (x) =Re \{ L^2 \cdots L^{l-2} \, \partial r \, ([L^{l-1},L^l]) (x)\}$$ and $$g (x) =Im \{ L^2 \cdots L^{l-2} \, \partial r \, ([L^{l-1},L^l]) (x)\}.$$ The condition $\ls \partial r (x_0)\neq 0$ implies  at least one of $Xf(x_0),$ $Xg(x_0),$ $Yf(x_0),$ and $Yg(x_0)$ cannot vanish. Without loss of generality, let us assume $Xf(x_0) \neq 0.$ We set $r_{p+2}(x)=f(x)$ and $L_{p+2}=L,$ the vector field from which we constructed the list $\ls.$ Inductively, assume that for integer $\nu-1$ with $p+2 \leq \nu-1<n,$ we have already constructed entries $c_1, \dots, c_{\nu-1};$ functions $r_1,$ $r_{p+2}, \dots, r_{\nu-1};$ and vector fields $L_2, \dots, L_{\nu-1}.$ Denote by $T^{\, (1,0)}_\nu$ the set of $(1,0)$ smooth vector fields $L$ such that $\partial \bar \partial r (L, \bar L_j)=0$ for $j=2, \dots, p+1$ and $L(r_k)=0$ for $k=1,p+2, p+3, \dots,\nu-1.$ Just as before, let $\T_\nu$ be the set of germs of sections of $T^{\, (1,0)}_\nu.$ Fix a vector field $L$ in $\T_\nu,$ and consider the list $\ls=\{L^1, \dots, L^l\}$ satisfying that each $L^i$ is one of the vector fields from the set $\{L_{p+2}, \bar L_{p+2}, \dots, L_{\nu-1}, \bar L_{\nu-1},L, \bar L\}.$ Let $l_i$ denote the total number of times both $L_i$ and its conjugate $\bar L_i$ occur in $\ls$ for $p+2 \leq i \leq \nu-1,$ and let $l_\nu$ denote the total number of times both $L$ and $\bar L$ occur in the list $\ls.$ We will only consider lists $\ls$ that satisfy the following two definitions:

\medskip
\newtheorem{ordlist}[wgt]{Definition}
\begin{ordlist}
\label{ordlistdef} A list $\ls= \{L^1, \dots, L^l\}$ is called ordered if
\begin{enumerate}
\item[(i)] $L^j=L$ or $L^j= \bar L$ for $1 \leq j \leq l_\nu$
\item[(ii)] $L^j=L_i$ or $L^j=\bar L_i$ for $1+\sum_{k=i+1}^\nu l_k \leq j \leq \sum_{k=i}^\nu l_k.$
\end{enumerate}
\end{ordlist}

\medskip
\newtheorem{admlist}[wgt]{Definition}
\begin{admlist}
\label{admlistdef} A list $\ls= \{L^1, \dots, L^l\}$ is called $j$-admissible if
\begin{enumerate}
\item[(i)] $l_j>0;$
\item[(ii)] $$\sum_{i=p+2}^{j-1} \frac{l_i}{c_i}<1,$$ where $\ct^{j-1}=(c_1, \dots, c_{j-1})$ is the $(j-1)^{th}$ commutator multitype.
\end{enumerate}
\end{admlist}

\smallskip\noindent See \cite{levidet} or \cite{catlinbdry} for more motivation regarding these definitions. The content of Catlin's Theorem 6.3 from p.552 of \cite{catlinbdry} works perfectly as a definition of a boundary system except for the assumption of pseudoconvexity, which is not necessary:

\medskip
\newtheorem{bdrysys}[wgt]{Definition}
\begin{bdrysys}
Let $\B_\nu (x_0) = \{r_1, r_{p+2}, \dots, r_\nu; L_2, \dots, L_\nu\}$ be a set of real-valued smooth functions and vector fields in a neighborhood of a point $x_0$ in $\C^n,$ where $p+1 \leq \nu \leq n.$ Assume that the following properties are satisfied:
\begin{enumerate}
\item Near $x_0,$ $r_1$ is the defining function of a smoothly bounded domain. If we set $a_{ij}(x) =\partial \bar \partial r_1(L_i, L_j)(x),$ then the $p \times p$ Hermitian matrix $(a_{ij}(x_0))$ for $2 \leq i,j \leq p+1$ is nonsingular. If either $i$ or $j$ is greater than $p+1,$ then $a_{ij}(x_0)=0.$
\item If $k<j,$ then $L_j r_k \equiv 0.$ Also, the vector fields $L_2, \dots, L_\nu$ are linearly independent.
\item Let $(c_1, \dots, c_\nu)$ be a given weight in $\Gamma_\nu$ with $c_1=1,$ $c_2= \cdots = c_{p+1}=2,$ and $c_i \geq 3$ for $p+2 \leq i \leq \nu.$ For every $j=p+2, \dots, \nu,$ there is a $j$-admissible ordered list $\ls_j= \{L^1, \dots, L^l\}$ with $L^1=L_j$ or $\bar L_j$ such that $\ls_j \partial r_1 (x_0) \neq 0.$ If $l^j_i$ equals the number of times $L_i$ and $\bar L_i$ occur in $\ls_j,$ then $\sum_{i=p+2}^j \frac{l^j_i}{c_i}=1.$ Moreover, if $\ls'_j$ is defined by $\{L^2, \dots, L^l\},$ then $r_j$ equals either $Re\{\ls'_j \partial r_1\}$ or $Im\{\ls'_j \partial r_1\}.$
\item Let $\ls$ be any ordered list. If $l_i$ equals the number of times $L_i$ and $\bar L_i$ occur in $\ls$ for $p+2 \leq i \leq \nu$ and if $\sum_{i=p+2}^\nu \frac{l_i}{c_i}<1,$ then $\ls \partial r_1 (x_0) =0.$
\end{enumerate}
Then under these assumptions, $\B_\nu (x_0) = \{r_1, r_{p+2}, \dots, r_\nu; L_2, \dots, L_\nu\}$ is a boundary system of rank $p$ and codimension $n-\nu$ about the point $x_0.$ The $\nu^{th}$ commutator multitype of the boundary system $\B_\nu$ at $x_0$ is the $\nu$-tuple $\ct^\nu=(c_1, \dots, c_\nu),$ and vector fields $L_2, \dots, L_\nu$ are called the special vector fields associated to the boundary system $\B_\nu.$
\label{bdrysysdef}
\end{bdrysys}

\smallskip\noindent {\bf Remarks:}

\smallskip\noindent (1) As explained in the paragraph preceding Definition~\ref{ordlistdef}, $c_{p+2}$ is a positive integer, and $\ls_{p+2}$ only consists of $L_{p+2}$ and $\bar L_{p+2}.$

\smallskip\noindent (2) $\ct^\nu$ is well-defined, i.e. independent of coordinate system and choices of $r_{p+2}, \dots, r_\nu, L_2, \dots, L_\nu.$ In other words, the lengths of the admissible lists stay the same, but we may choose the vector fields $L_2, \dots, L_\nu$ slightly differently, which in turns would yield different functions $r_{p+2}, \dots, r_\nu$ in the boundary system.

\medskip\noindent We will make use of the freedom hinted at in the last remark in order to put the Levi form around $x_0$ in as close to a diagonalized form as this setting allows. The following lemma is an argument communicated to the author by David Catlin, although in a less transparent way it already appears on pp.539-540 of \cite{catlinbdry} and is listed among the properties of a boundary system on p.552 of the same paper:

\medskip
\newtheorem{nearlydiag}[wgt]{Lemma}
\begin{nearlydiag}
If $\B_\nu (x_0) = \{r_1, r_{p+2}, \dots, r_\nu; L_2, \dots, L_\nu\}$ is a boundary system of rank $p$ and codimension $n-\nu$ about the point $x_0 \in b \Omega,$ then there exists another boundary system $\tilde \B_\nu (x_0) = \{r_1, \tilde r_{p+2}, \dots, \tilde r_\nu; \tilde L_2, \dots, \tilde L_\nu\}$ at $x_0$ that is also of rank $p$ and codimension $n-\nu$ and furthermore satisfies that $\tilde a_{ij}(x) =\partial \bar \partial r_1 (\tilde L_i, \tilde L_j)(x)\equiv 0$ in a neighborhood of $x_0$ whenever $2 \leq i \leq p+1$ and $j \geq p+2.$
\label{nearlydiaglemma}
\end{nearlydiag}

\smallskip\noindent {\bf Proof:} The argument proceeds in two steps.

\smallskip\noindent {\bf Step 1:} Since $\B_\nu (x_0) = \{r_1, r_{p+2}, \dots, r_\nu; L_2, \dots, L_\nu\}$ is a boundary system of rank $p$ at $x_0,$ the Levi form has rank $p$ at $x_0,$ so there exist vector fields $\tilde L_2, \dots, \tilde L_{p+1}$ in a neighborhood $U$ of $x_0$ in $b \Omega$ such that the $p \times p$ Hermitian matrix $(\tilde a_{ij}(x))$ for $\tilde a_{ij}(x) =\partial \bar \partial r_1 (\tilde L_i, \tilde L_j)(x)$ and $2 \leq i,j \leq p+1$ is nonsingular at every $x \in U$ and furthermore equals the identity matrix $I_p$ at $x_0$ itself.

\smallskip\noindent {\bf Step 2:} Complete $\tilde L_2, \dots, \tilde L_{p+1}$ to a basis $\tilde L_2, \dots, \tilde L_{p+1}, L'_{p+2} \dots, L'_n$ of $T^{\, (1,0)} (U \cap b \Omega).$ We claim that for every $k, l$ with $p+2 \leq k \leq n$ and $2 \leq l \leq p+1,$ there exist smooth functions $c_{k,l} \in \smooth (U)$ such that the vector fields $\tilde L_k = L'_k + \sum_{l=2}^{p+1} \, c_{k,l} \, \tilde L_l$ satisfy $\tilde a_{ik}(x) =\partial \bar \partial r_1 (\tilde L_i, \tilde L_k)(x)\equiv 0$ for $2 \leq i \leq p+1$ on the neighborhood $U$ of $x_0$ from Step 1 up to perhaps a shrinking.

\smallskip\noindent {\bf Proof of claim:} $\tilde a_{ik}(x) =\partial \bar \partial r_1 (\tilde L_i, \tilde L_k)(x)\equiv 0$ for all $x \in U$ is equivalent to $$\partial \bar \partial r_1 \left(\tilde L_i, L'_k + \sum_{l=2}^{p+1} \, c_{k,l} \, \tilde L_l\right)(x)\equiv 0,$$ which is in turn equivalent to
\begin{equation}
 \partial \bar \partial r_1 \left(\tilde L_i, \sum_{l=2}^{p+1} \, c_{k,l} \, \tilde L_l\right)(x)=- \partial \bar \partial r_1 (\tilde L_i,  L'_k)(x). \label{leviformeq}
\end{equation}
Set $b_{i,k} = - \partial \bar \partial r_1 (\tilde L_i,  L'_k)(x)$ for $2 \leq i \leq p+1$ and $p+2 \leq k \leq n.$ Obviously, these $b_{i,k}$'s are smooth functions in $\smooth(U).$ By the linearity of the pairing of forms and vector fields that defines the Levi form, for every $k$ such that $p+2 \leq k \leq n,$ equation ~\eqref{leviformeq} can be rewritten as the linear system of equations
\begin{equation}
\left(\partial \bar \partial r_1(\tilde L_i,  \tilde L_l)\right)
\begin{pmatrix}
c_{k,2}\\
\vdots\\
c_{k,p+1}
\end{pmatrix}
=
\begin{pmatrix}
b_{2,k}\\
\vdots\\
b_{p+1,k}
\end{pmatrix}, \label{levisys}
\end{equation}
where the $p \times p$ Hermitian matrix $\left(\partial \bar \partial r_1(\tilde L_i,  \tilde L_l)\right)$ is the identity matrix $I_p$ at $x_0$ by Step 1 and a small perturbation of the identity in a neighborhood of $x_0.$ Shrink the neighborhood $U$ from Step 1 as necessary to ensure that by Cramer's rule, the linear system of equations in \eqref{levisys} has smooth solutions $c_{k,2}, \dots, c_{k,p+1} \in \smooth(U)$ for every $k$ such that $p+2 \leq k \leq n.$ The claim is proven. 

Note that while the defining function $r_1$ does not change in the new boundary system with respect to vector fields $\tilde L_2, \dots, \tilde L_\nu,$ we will be obtaining different functions $\tilde r_{p+2}, \dots, \tilde r_\nu.$ Additionally, a quick glance at the definition of the tangent bundle $T^{\, (1,0)}_{p+2}$ shows that however we choose $L'_{p+2} \dots, L'_n \in T^{\, (1,0)}_{p+2}$ to complete $\tilde L_2, \dots, \tilde L_{p+1}$ to a basis of $T^{\, (1,0)} (U \cap b \Omega),$ the new vector fields given by $\tilde L_k = L'_k + \sum_{l=2}^{p+1} \, c_{k,l} \, \tilde L_l$ also satisfy $\tilde L_k \in T^{\, (1,0)}_{p+2}$ for every $k$ such that $p+2 \leq k \leq n.$ \qed

\medskip\noindent From now on, we can assume our boundary system $\B_\nu (x_0) = \{r_1, r_{p+2}, \dots, r_\nu; L_2, \dots, L_\nu\}$ satisfies the conclusion of Lemma~\ref{nearlydiaglemma} whenever it is useful to do so.

We summarize in the next theorem two of the most important properties of $\ct^\nu,$ which are contained in Proposition 2.1 on p.536 and Theorem 2.2 on p.538 of \cite{catlinbdry}:

\medskip
\newtheorem{cmtype}[wgt]{Theorem}
\begin{cmtype}
\label{cmtypethm} Let the domain $\Omega= \{ z \in \C^n \: \big| \: r(z)<0\}$ be pseudoconvex in a neighborhood of a point $x_0 \in b \Omega.$ The $\nu^{th}$ commutator multitype $\ct^\nu=(c_1, \dots, c_\nu)$ of the boundary system $\B_\nu$ satisfies the following two properties:
\begin{enumerate}
\item[(i)] $\ct^\nu$ is upper semi-continuous with respect to the lexicographic ordering, i.e. there exists a neighborhood $U \ni x_0$ such that for all $x \in U \cap b \Omega,$ $\ct^\nu (x) \leq \ct^\nu (x_0).$
\item[(ii)] $\ct^\nu(x_0)=\mt^\nu(x_0),$ where $\mt^\nu=(m_1, \dots, m_\nu)$ consists of the first $\nu$ entries of the multitype $\mt=(m_1, \dots, m_n).$
\end{enumerate}
\end{cmtype}

\smallskip\noindent {\bf Remark:} Pseudoconvexity is not necessary for part (i) but is essential for part (ii) of this theorem.

We will now state only part of the main theorem on p.531 of \cite{catlinbdry} containing the properties of the multitype $\mt(x_0):$

\medskip
\newtheorem{mltthm}[wgt]{Theorem}
\begin{mltthm}
\label{mlttheorem} Let $\Omega \subset \C^n$ be a pseudoconvex domain with smooth boundary. Let $x_0 \in b \Omega.$ The multitype $\mt(x_0)$ has the following properties:
\begin{enumerate}
\item If $\mt(x_0)=(m_1, \dots, m_n),$ then there exist coordinates $(z_1, \dots, z_n)$ around $x_0$ such that $x_0$ is mapped to the origin and if $\sum_{i=1}^n \: \frac{\alpha_i+ \beta_i}{m_i}<1,$ then $D^\alpha \bar D^\beta r(0)=0.$ If one of the entries $m_i = + \infty$ for some $1 \leq i \leq n,$ then these coordinates should be interpreted in the sense of formal power series.
\item If $\mt(x_0)=(m_1, \dots, m_n),$ then for each $q=1, \dots, n,$ $$m_{n+1-q} \leq \Delta_q (b \Omega, x_0),$$ where $\Delta_q (b \Omega, x_0)$ is the D'Angelo $q$-type of the point $x_0.$
\end{enumerate}
\end{mltthm}

\medskip In Section~\ref{effcomp}, we will need to truncate the defining function of the domain in a way that retains finite D'Angelo type but may lose pseudoconvexity so that we can understand more about what types of derivatives of the Levi determinant are non-zero. Therefore, some machinery from \cite{catlinbdry} that works in the absence of pseudoconvexity will be very useful. We recall it here starting with a definition that introduces a measurement of the vanishing order of a function with respect to a given weight:

\medskip
\newtheorem{disfct}[wgt]{Definition}
\begin{disfct}
\label{disfctdef} Let the weight $\Lambda=(\lambda_1, \dots, \lambda_n) \in \Gamma_n$ be given. We will denote by $\mt(t; \Lambda)$ the set of germs of smooth functions $f$ defined near the origin such that  $$D^\alpha \bar D^\beta f(0)=0 \:\: \text{whenever} \:\:\sum_{i=1}^n \: \frac{\alpha_i+ \beta_i}{\lambda_i}<t.$$
\end{disfct} 

\medskip We will now state Catlin's Proposition 3.6 from page 542 of \cite{catlinbdry}. Given a boundary system at $x_0,$ $$\B_\nu(x_0)=\{r_1,r_{p+2}, \dots, r_\nu; L_2, \dots, L_\nu\},$$ this result shows in which sets $\mt(t;\Lambda)$ we can place the additional functions $r_{p+2}, \dots, r_\nu$ obtained by differentiating $r$ and relates this information to the $\nu^{th}$ commutator multitype $\ct^\nu (x_0).$ Without loss of generality, we can translate the point $x_0$ to the origin in $\C^n.$

\medskip
\newtheorem{placerj}[wgt]{Proposition}
\begin{placerj}
\label{placerjprop} Let $\B_\nu(0)=\{r_1,r_{p+2}, \dots, r_\nu; L_2, \dots, L_\nu\}$ be a boundary system of rank $p$ and codimension $n-\nu$ about the origin in $\C^n.$ Suppose that the $\nu^{th}$ commutator multitype of $\B_\nu(0)$ at the origin $\ct^\nu(0)=(\lambda_1, \dots, \lambda_\nu)$ and let $\Lambda=(\lambda_1, \dots, \lambda_n)$ be a weight in $\Gamma_n$ that agrees with $\ct^\nu(0)$ up to the $\nu^{th}$ entry and also satisfies that $\lambda_1=1,$ $\lambda_2= \cdots = \lambda_{p+1}=2,$ and $\lambda_j \geq 3$ for $j \geq p+2.$ If $r_1 \in \mt(1; \Lambda),$ then $r_k \in \mt\left(\frac{1}{\lambda_k};\Lambda\right)$ for all $k=p+2, \dots, \nu.$ Moreover, if $c_{\nu+1}$ denotes the $(\nu+1)^{th}$ entry of $\ct^{\nu+1}(0),$ then $c_{\nu+1} \geq \lambda_{\nu+1}.$
\end{placerj}

\smallskip\noindent As Proposition~\ref{placerjprop} shows, the $\nu$-tuple $\ct^\nu(0)$ is associated to the boundary system $\B_\nu(0),$ but the space $\mt(t;\Lambda)$  requires a weight $\Lambda \in \Gamma_n,$ which is an $n$-tuple. We thus need to manufacture a weight starting with an $\nu$-tuple. The most natural way to do so is contained in the next definition:

\medskip
\newtheorem{elongate}[wgt]{Definition}
\begin{elongate}
\label{elongatedef} Let $\Gamma_{\, n, \nu+1}$ be the set of weights $(\lambda_1, \dots, \lambda_n)$ in $\Gamma_n$ such that $\lambda_{\nu+1} = \cdots = \lambda_n,$ i.e. all entries from the $(\nu+1)^{th}$ entry forward coincide.
\end{elongate}

\smallskip Definition~\ref{elongatedef} will become relevant in Section~\ref{effcomp}, if it turns out that $\mt(x_0)$ has some infinite entries as it could be the case if $\Delta_q (b \Omega, x_0)<\infty,$ but $q>1.$ 

We can finally prove now the lemma that will be used in Section~\ref{effcomp} to identify the non-zero derivatives of the Levi determinant from which condition (a) in the Conjectured Nullstellensatz~\ref{conjecturednullstellensatz} can be derived. The proof of this lemma is part of the proof of Catlin's Proposition 3.2 on p.539 of \cite{catlinbdry}, which we simply employ here for a different purpose compared to Catlin:

\medskip
\newtheorem{newrj}[wgt]{Lemma}
\begin{newrj}
\label{newrjlemma} Let $\B_\nu(0)=\{r_1,r_{p+2}, \dots, r_\nu; L_2, \dots, L_\nu\}$ be a boundary system of rank $p$ and codimension $n-\nu$ about the origin in $\C^n.$ Suppose that the $\nu^{th}$ commutator multitype of $\B_\nu(0)$ at the origin $\ct^\nu(0)=(\lambda_1, \dots, \lambda_\nu)$ and let $\Lambda=(\lambda_1, \dots, \lambda_n)$ be a weight in $\Gamma_n$ that agrees with $\ct^\nu(0)$ up to the $\nu^{th}$ entry and also satisfies that $\lambda_1=1,$ $\lambda_2= \cdots = \lambda_{p+1}=2,$ and $\lambda_j \geq 3$ for $j \geq p+2.$ If $r_k \in \mt\left(\frac{1}{\lambda_k};\Lambda\right)$ for all $k=1, p+2, \dots, \nu,$ then there exists a coordinate change $w=\psi(z)$ in a neighborhood of the origin in $\C^n$ such that $\rho_k=(\psi^{-1})^* r_k,$ the pullbacks of the functions $r_k$ under this coordinate change for $k=1, p+2, \dots, \nu,$ satisfy $\rho_k \in \mt\left(\frac{1}{\lambda_k};\Lambda\right)$ and for $p+2\leq k \leq \nu,$ $$\rho_k(w) = 2 Re\{w_k\} + O(|w|^2).$$ Furthermore, if $\tilde \B_\nu(0)=\{\rho_1,\rho_{p+2}, \dots, \rho_\nu; \tilde L_2, \dots, \tilde L_\nu\}$ is the boundary system of rank $p$ and codimension $n-\nu$ about the origin in $\C^n$ in the new coordinates corresponding to $\B_\nu(0),$ then for each $k$ such that $p+2 \leq k \leq \nu,$ the coefficient of $\frac{\partial}{\partial w_k}$ in the vector field $\tilde L_k$ is non-zero at the origin.
\end{newrj}

\smallskip\noindent {\bf Proof:} Parts (2) and (3) of Definition~\ref{bdrysysdef} imply that the $(\nu-p-1) \times (\nu-p-1)$ matrix $(L_i r_j(0)),$ where $p+2 \leq i,j \leq \nu,$ is upper triangular and that its diagonal entries are all non-zero. Therefore, we conclude that at the origin the differentials $$\partial r_k (0) = \sum_{j=1}^n \, a_j^k \, d z_j$$ are linearly independent for $k=p+2, \dots, \nu.$ Furthermore, the assumption $r_k \in \mt\left(\frac{1}{\lambda_k};\Lambda\right)$ for all $k=p+2, \dots, \nu$ implies $a_j^k=0$ if $\lambda_j > \lambda_k.$ We now choose additional complex numbers $a_j^k$ for $\nu+1 \leq k \leq n$ and $1 \leq j \leq n$ such that $a_j^k=0$ if $\lambda_j > \lambda_k$ and the matrix $(a_j^k)$ has rank $n-p-1$ for $p+2 \leq k \leq n$ and $1 \leq j \leq n.$ To obtain a full rank matrix $(a_j^k)$ for $1 \leq j,k \leq n,$ we set $a_j^k=\delta_{jk}$ whenever $1 \leq k \leq p+1$ and $1 \leq j \leq n.$ We use this matrix to define a linear change of variables around the origin in $\C^n$ as follows: $w= \psi(z),$ where for every $k=1, \dots, n,$ $$\psi_k(z)= \sum_{j=1}^n \, a_j^k z_j.$$ Consider now $\rho_k=(\psi^{-1})^* r_k,$ the pullbacks of the functions $r_k$ under this coordinate change for $k=1, p+2, \dots, \nu.$ The condition we imposed that $a_j^k=0$ if $\lambda_j > \lambda_k$ ensures $\rho_k \in \mt\left(\frac{1}{\lambda_k};\Lambda\right)$ for $k=1, p+2, \dots, \nu.$ The inverse $\psi^{-1}$ is a linear map represented by the matrix that is the inverse of $(a_j^k),$ whose rows $p+2$ through $\nu$ were exactly the differentials $\partial r_{p+2}(0), \dots, \partial r_\nu (0).$ Therefore, $\partial \rho_k (0) = d w_k$ and $\rho_k(w) = 2 Re\{w_k\} + O(|w|^2)$ in a neighborhood of the origin for $k=p+2, \dots, \nu.$ To the original boundary system $\B_\nu(0)=\{r_1,r_{p+2}, \dots, r_\nu; L_2, \dots, L_\nu\},$ there corresponds a boundary system $\tilde \B_\nu(0)=\{\rho_1,\rho_{p+2}, \dots, \rho_\nu; \tilde L_2, \dots, \tilde L_\nu\}$ also of rank $p$ and codimension $n-\nu$ about the origin in $\C^n.$ Therefore, $\tilde L_k \rho_k (0) \neq 0$ for every $k=p+2, \dots, \nu,$ which given the form of $\rho_k$ around the origin means the coefficient of $\frac{\partial}{\partial w_k}$ in the vector field $\tilde L_k$ must be non-zero at the origin. \qed

\medskip\noindent We conclude this section with a lemma that follows easily from Catlin's construction of a boundary system but is essential for the proof of Theorem~\ref{maintheorem}:

\medskip
\newtheorem{point3}[wgt]{Lemma}
\begin{point3}
\label{point3lemma} Let the domain $\Omega= \{ z \in \C^n \: \big| \: r(z)<0\}$ be smooth in a neighborhood of a point $x_0 \in b \Omega,$ and let $$\B_{n+1-q} (x_0) = \{r_1, r_{p+2}, \dots, r_{n+1-q}; L_2, \dots, L_{n+1-q}\}$$ be a boundary system of rank $p$ and codimension $n-(n+1-q)=q-1$ at $x_0,$ then $$\partial r \wedge \bar \partial r \wedge (\partial \bar \partial r)^p \wedge \partial r_{p+2} \wedge \cdots \wedge \partial r_{n+1-q} (x_0) \neq 0.$$
\end{point3}

\smallskip\noindent {\bf Remark:} If we show $ r_{p+2}, \dots, r_{n+1-q}$ are multipliers in the ideal $I^q_k(x_0)$ at some step $k$ of the Kohn algorithm, then $$\text{coeff}\{\partial r \wedge \bar \partial r \wedge (\partial \bar \partial r)^p \wedge \partial r_{p+2} \wedge \cdots \wedge \partial r_{n+1-q} \} \in A^q_k (x_0)$$ because $n+1-q-(p+1)+p=n-q.$ Thus $\text{coeff}\{\partial r \wedge \bar \partial r \wedge (\partial \bar \partial r)^p \wedge \partial r_{p+2} \wedge \cdots \wedge \partial r_{n+1-q} \} $ would be a unit in the ideal $I^q_{k+1} (x_0),$ and the Kohn algorithm would terminate at step $k+1.$

\smallskip\noindent {\bf Proof:} From parts (2) and (3) of Definition~\ref{bdrysysdef}, for $i=1, p+2, p+3, \dots, n+1-q$ and $j=2, \dots, n+1-q,$ we have $$L_j (r_i)\begin{cases} =0 & \: \text{if} \: j>i\\  \neq 0 & \: \text{if} \: j=i\\  \text{no information} & \: \text{if} \: j<i \end{cases}.$$ The vector fields $L_2, \dots, L_{n+1-q}$ are linearly independent and belong to $T^{\, (1,0)} (b \Omega \cap U)$ for $U \ni x_0$ an open set around $x_0.$ We now complete these to a basis of $T^{\, (1,0)} (b \Omega)$ in accordance with the setup on page~\pageref{prenbd}. Since the imaginary part of $L_1$ is $T$ and its real part is the normal direction to $b \Omega$ as shown on page~\pageref{prenbd}, it follows that at $x_0$ the wedge product $\partial r \wedge \bar \partial r \wedge (\partial \bar \partial r)^p \wedge \partial r_{p+2} \wedge \cdots \wedge \partial r_{n+1-q}$ is given by all the $(n-q) \times (n-q)$ minors of the $(n-1) \times (n-q)$ matrix:
\begin{equation*}
\begin{pmatrix}
A_p & 0 \\ \star & B
\end{pmatrix},
\end{equation*}
where $A_p$ is a $p \times p$ nonsingular matrix coming from the fact that the Levi form has rank $p$ at $x_0,$ $0$ is a $(n-1-p) \times p$ matrix of all zero entries, $\star$ is a $p \times (n-q-p)$ matrix for which we have no information, and $B$ is the following lower triangular $(n-1-p) \times (n-q-p)$ matrix:
\begin{equation*}
\begin{pmatrix}
L_{p+2}(r_{p+2}) & 0 & \cdots &0 & 0 & 0 & \cdots & 0 \\ \ast & L_{p+3}(r_{p+3}) & \cdots &0 & 0  & 0 & \cdots & 0\\ \vdots & \vdots & \vdots & \vdots & \vdots & \vdots & \vdots & \vdots \\ \ast & \ast & \cdots & L_{n-q}(r_{n-q})  &0  & 0 & \cdots & 0 \\ \ast & \ast & \cdots & \ast & L_{n+1-q}(r_{n+1-q})  & 0 & \cdots & 0
\end{pmatrix}
\end{equation*}
Notation $\ast$ denotes an entry for which we have no information. Note that the right side block of zero entries occurs only if $q >1.$ Given the location of the zero entries, it is clear that the wedge product $$\big(\partial r \wedge \bar \partial r \wedge (\partial \bar \partial r)^p \wedge \partial r_{p+2} \wedge \cdots \wedge \partial r_{n+1-q}\big) (x_0) =( \det A_p) \, L_{p+2}(r_{p+2}) \cdots  L_{n+1-q}(r_{n+1-q}) \neq 0$$ by construction. \qed

\section{Tougeron-Whitney elements and the stratification theorem}
\label{stratsect}

The aim of this section is to define Tougeron-Whitney elements, relate Tougeron-Whitney elements to Catlin's functions $r_1, r_{p+2}, \dots, r_{n+1-q}$ in a boundary system $\B_{n+1-q},$ and then prove a stratification theorem that highlights the importance of the observation contained in Lemma~\ref{point3lemma}. The author also utilized ideas related to the stratification induced by the Catlin multitype in order to investigate the behavior of the Levi determinant in \cite{levidet}.

Whitney is said to have first used the objects we will be defining shortly, and they subsequently appear in work by Malgrange and others. In particular, in Proposition $4.6$ of subsection $V.4$ of \cite{tougeron}, Tougeron employed these objects to simplify the proof of Thom's result in \cite{thom} that a variety corresponding to a {\L}ojasiewicz ideal of smooth functions has an open and dense set of smooth points. Tougeron's use is closest to the situation at hand, so we have decided to call these objects Tougeron-Whitney elements. Since these elements are $\R$-valued smooth functions, if $k=(k_1, \dots, k_n) \in \N^n$ is a multi-index, then $D^k$ means the differentiation
$\frac{\partial}{\partial x_1^{k_1}}\cdots
\frac{\partial}{\partial x_n^{k_n}}$ in this context.

\smallskip
\newtheorem{tougeronelem}{Definition}[section]
\begin{tougeronelem}
Consider $f \in \smooth (\R^n)$ and $\V=\{x \in \R^n \: | \:
f(x)=0\}= \V \big( (f) \big).$ If $x_0 \in \V$ and there exist an
open set $U\subset \R^n,$ $x_0 \in U,$ and a multi-index $k \in
\N^n$ with $|k|=d\geq 1$ such that $D^k f(x_0) \neq 0$ but $D^h
f(x)=0$ for all $h \in \N^n$ such that $|h|<|k|=d$ and all $x \in
U \cap \V,$ then we call $g(x)=D^{k'} f(x)$ a Tougeron-Whitney element
corresponding to $f(x)$ provided that the multi-index $k' \in
\N^n$ arises from the multi-index $k$ by splitting off a factor of $\frac{\partial}{\partial x_j},$ i.e. $\frac{\partial}{\partial x_j}
D^{k'} f(x_0)=D^k f(x_0)\neq 0$ for some $1 \leq j \leq n.$ \label{tougeronelemdef}
\end{tougeronelem}

\smallskip\noindent {\bf Remark:} By construction, the gradient of the
Tougeron-Whitney element $g(x)$ satisfies that $\nabla g(x_0) \neq 0,$ so
there exists a perhaps smaller open set $\tilde U \subset U$ with
$x_0 \in \tilde U$ such that $ \V \big( (g) \big) \cap \tilde U$
is a hypersurface, where $ \V \big( (g) \big)$ is the vanishing
set of $g(x).$

Tougeron's setup in \cite{tougeron} was the author's motivation for looking at Catlin's boundary system construction in \cite{catlinbdry} in the hope of finding distinguished elements with nonzero gradients that vanished on the varieties corresponding to the ideals $I^q_k (x_0)$ in the Kohn algorithm. We shall thus call functions $r_1, r_{p+2}, \dots, r_{n+1-q}$ in the boundary system $\B_{n+1-q}(x_0)$ Tougeron-Whitney elements for the Kohn algorithm. A consequence of Lemma~\ref{point3lemma} in the previous section is that functions $r_1, r_{p+2}, \dots, r_{n+1-q}$ have linearly independent, nonzero gradients. We now have to show that indeed $r_1, r_{p+2}, \dots, r_{n+1-q}$ vanish on the variety corresponding to the first ideal of multipliers $I^q_1(U),$ provided we subtract certain level sets of the $(n+1-q)^{th}$ commutator multitype $\ct^{n+1-q}$ and $U$ is an appropriately small neighborhood around $x_0.$ This result is the stratification theorem mentioned at the beginning of this section.

To use the stratification given by Catlin's multitype, just like in \cite{levidet}, we need a beefed-up version of Proposition 2.1 on p.536 of \cite{catlinbdry}

\medskip
\newtheorem{stratcatlin}[tougeronelem]{Proposition}
\begin{stratcatlin}
\label{stratcatlinprop} Let $x_0 \in b \Omega$ be such that the Levi form has rank $p$ at $x_0.$ For $p+2 \leq \nu \leq n,$ let $\B_\nu$ be a boundary system of rank $p$ and codimension $n- \nu$ at $x_0.$ There exists a neighborhood $U$ of $x_0$ such that all the following conditions hold on its closure $\overline{U}$:
\begin{enumerate}
\item[(i)] For all $x \in \overline{U} \cap b \Omega,$ $\ct^\nu (x) \leq \ct^\nu (x_0),$ where $\ct^\nu=(c_1, \dots, c_\nu)$ is the $\nu^{th}$ commutator multitype;
\item[(ii)] $$M^\nu=\{ x \in \overline{U} \cap b \Omega \: \big| \: r_j(x)=0, \: j=1,p+2, \dots, \nu\}$$ is a submanifold of $\overline{U} \cap b \Omega$ of holomorphic dimension $n-\nu$ in the sense of Catlin;
\item[(iii)] The level set of the commutator multitype at $x_0$ satisfies that $$\{ x \in \overline{U} \cap b \Omega \: \big| \: \ct^\nu(x)=\ct^\nu(x_0)\} \subset M^\nu;$$
\item[(iv)] For all $x \in \overline{U} \cap b \Omega,$ the Levi form has rank at least $p$ at $x;$
\item[(v)] For all $x \in \overline{U} \cap b \Omega,$ $\ls_j \partial r_1(x) \neq 0$ for all $j=p+2, \dots, \nu,$ where $\ls_{p+2}, \dots, \ls_\nu$ are the $\nu$-admissible, ordered lists used in defining the boundary system $\B_\nu.$
\end{enumerate}
\end{stratcatlin}

\smallskip Now let $x_0 \in b\Omega$ be a point of finite D'Angelo $q$-type, and assume $\Omega$ is pseudoconvex. By D'Angelo's result, Theorem~\ref{deltaqopenthm}, there exists a neighborhood $U$ of $x_0$ in $b \Omega$ so that for every $x \in U,$ the D'Angelo $q$-type is finite at $x.$ We can shrink $U$ if necessary to ensure Proposition~\ref{stratcatlinprop} also holds on $U.$ By remark 1.2 on p.532 of \cite{catlinbdry},  the $(n+1-q)^{th}$ commutator multitype $\ct^{n+1-q}$ can assume only finitely many values $\ct^{n+1-q}_1, \dots, \ct^{n+1-q}_N$ at all points of $U \cap b \Omega,$ where $\ct^{n+1-q}_1< \ct^{n+1-q}_2 < \cdots < \ct^{n+1-q}_N$ and $N$ is some positive natural number. We will give an effective upper bound for $N$ in terms of the D'Angelo $q$-type, $n,$ and $q$ later on in the paper. Let $$S_j= \{ x \in U \cap b \Omega \: \big| \: \ct^{n+1-q} (x)=\ct^{n+1-q}_j\}$$ be the level sets of the $(n+1-q)^{th}$ commutator multitype for $1 \leq j \leq N.$ We are working here with the open set $U$ rather than its closure, unlike in Proposition~\ref{stratcatlinprop}. We now have the stratification $$U \cap b \Omega = \bigcup_{j=1}^N \, S_j$$ since $S_i \cap S_j = \emptyset$ for $i \neq j.$

\smallskip Let us recall Lemmas 4.9 and 4.10 from \cite{levidet} by combining them into one statement:

\medskip
\newtheorem{lowopen}[tougeronelem]{Lemma}
\begin{lowopen}
\label{lowopenlemma} Let $\Omega \subset \C^n$ be a smooth pseudoconvex domain, and let $x_0 \in b \Omega$ be a boundary point of finite D'Angelo $q$-type. Let $U$ be a neighborhood of $x_0$ such that on $U \cap b \Omega,$ $b \Omega$ has finite D'Angelo $q$-type everywhere and the $(n+1-q)^{th}$ commutator multitype $\ct^{n+1-q}$ takes only finitely many values $\ct^{n+1-q}_1< \dots< \ct^{n+1-q}_N$ for some natural number $N \geq 1.$ The lowest $(n+1-q)^{th}$ commutator multitype $\ct^{n+1-q}_1=(1,2, \dots,2)$ and its level set $S_1$ is open in $b \Omega.$
\end{lowopen}

\smallskip
\newtheorem{firststep}[tougeronelem]{Corollary}
\begin{firststep}
\label{firststepcor} Let $\Omega \subset \C^n$ be a smooth pseudoconvex domain, and let $x_0 \in b \Omega$ be a boundary point of finite D'Angelo $q$-type. There exists a neighborhood $U$ of $x_0$ such that the Kohn algorithm terminates at step 1 densely in $U \cap b \Omega$ in the induced topology of $b \Omega.$
\end{firststep}

\smallskip\noindent {\bf Proof:} By Lemma~\ref{lowopenlemma}, $\ct^{n+1-q}_1=(1,2, \dots,2)$ and its level set $S_1$ is open. Therefore, $\partial r \wedge \bar \partial r \wedge (\partial \bar \partial r)^{n-q} \neq 0$ densely in $U\cap b \Omega,$ but $\text{coeff}\{\partial r \wedge \bar \partial r \wedge (\partial \bar \partial r)^{n-q}\} \in I^q_1,$ the first ideal of multipliers in the Kohn algorithm. Thus $I^q_1 (x) = \smooth(x)$ for a dense set in $U \cap b \Omega,$ i.e. the Kohn algorithm terminates at step 1 at each of the points in this dense set. \qed

\medskip\noindent We are finally ready to state and prove the most important result in this section, the stratification theorem:

\medskip
\newtheorem{ourstrat}[tougeronelem]{Stratification Theorem}
\begin{ourstrat}
\label{ourstratthm} Let $\Omega \subset \C^n$ be a smooth pseudoconvex domain, and let $x_0 \in b \Omega$ be a boundary point of finite D'Angelo $q$-type. Let  $$\B_{n+1-q}(x_0)=\{r_1,r_{p+2}, \dots, r_{n+1-q}; L_2, \dots, L_{n+1-q}\}$$ be the boundary system at $x_0.$ There exists a neighborhood $U$ of $x_0$ such that $$r_1,r_{p+2}, \dots, r_{n+1-q} \in \I \left(\V(I^q_1(U))-\bigcup_{j=2}^{N-1} S_j\right),$$ i.e. the functions $r_1,r_{p+2}, \dots, r_{n+1-q}$ vanish on the zero set of the first ideal of multipliers $I^q_1(U)$ of the Kohn algorithm on the neighborhood $U$ after we remove from the zero set all the level sets of the commutator multitype except for the lowest one and the highest one, which is the one at the point $x_0$ itself.
\end{ourstrat}

\smallskip\noindent {\bf Remark:} By the upper semi-continuity of the commutator multitype, the level set $S_N$ corresponding to the top commutator multitype in $U$ is a closed set in $b \Omega.$ When working over $\smooth(U),$ any closed set is a variety, which implies that $$\V(I^q_1(U))-\bigcup_{j=2}^{N-1} S_j= \V(I^q_1(U)) \cap S_N$$ is a variety as well, and therefore the notation $$r_1,r_{p+2}, \dots, r_{n+1-q} \in \I \left(\V(I^q_1(U))-\bigcup_{j=2}^{N-1} S_j\right)$$ makes sense.
 Assume now that $\Omega$ is real-analytic. It will be easy to see in the proof of Theorem~\ref{maintheorem} that in fact $$\V(I^q_1(U))-\bigcup_{j=2}^{N-1} S_j= \V(I^q_1(U)) \cap S_N$$ is a real-analytic variety. Furthermore, since the domain is real-analytic, both generators of $I^q_1$ are real-analytic and so are $r_1,r_{p+2}, \dots, r_{n+1-q}.$ It follows that the notation $$r_1,r_{p+2}, \dots, r_{n+1-q} \in \I \left(\V(I^q_1(U))-\bigcup_{j=2}^{N-1} S_j\right)$$ makes sense for the ideal of the variety $\V(I^q_1(U))-\bigcup_{j=2}^{N-1} S_j$ in $\ra (U)$ and even more correctly, $$r_1,r_{p+2}, \dots, r_{n+1-q} \in \I \left(\V(\tilde I^q_1(U))-\bigcup_{j=2}^{N-1} S_j\right).$$

\smallskip\noindent {\bf Proof:} We distinguish two cases:

\smallskip\noindent {\bf Case 1:} The Levi form has rank $n-q$ at $x_0.$ In this case, the commutator multitype $\ct^{n+1-q}(x_0)=(1, 2, \dots, 2)$ and the boundary system does not contain any other $r_j$'s besides the defining function $r_1:$ $$\B_{n+1-q}(x_0)=\{r_1; L_2, \dots, L_{n+1-q}\}.$$ We choose $U$ to be the neighborhood guaranteed by Lemma~\ref{lowopenlemma}. Since $r_1=r,$ by the very definition of $I^q_1,$ $r_1 \in I^q_1(U).$ No level sets of $\ct^{n+1-q}$ need to be removed because by the upper semi-continuity of $\ct^{n+1-q},$ it stays the same on all of $U,$ and it is the lowest one.

\smallskip\noindent {\bf Case 2:} The Levi form has rank $p$ with $p < n-q$ at $x_0.$ This means $\ct^{n+1-q}(x_0)>(1, 2, \dots, 2).$ Let $U$ be a neighborhood of $x_0$ such that Proposition~\ref{stratcatlinprop} holds and the D'Angelo $q$-type is finite for all $x \in U \cap b \Omega.$ $$U \cap b \Omega = \bigcup_{j=1}^N \, S_j,$$ where $\ct^{n+1-q}_1< \ct^{n+1-q}_2 < \cdots < \ct^{n+1-q}_N$ and $$S_j= \{ x \in U \cap b \Omega \: \big| \: \ct^{n+1-q} (x)=\ct^{n+1-q}_j\}$$ are the level sets of the $(n+1-q)^{th}$ commutator multitype for $1 \leq j \leq N.$ By part (i) of Proposition~\ref{stratcatlinprop}, $x_0 \in S_N.$ By part (iii) of Proposition~\ref{stratcatlinprop}, $S_N \subset M^{n+1-q},$ where $$M^{n+1-q}=\{ x \in U \cap b \Omega \: \big| \: r_j(x)=0, \: j=1,p+2, \dots, {n+1-q}\}.$$ When $x \in S_1,$ $\text{coeff}\{\partial r \wedge \bar \partial r \wedge (\partial \bar \partial r)^{n-q}\} \neq 0,$ so $S_1 \cap \V(I^q_1(U))= \emptyset.$ For all $2 \leq j \leq N$ and all $x \in S_j,$ $\text{coeff}\{\partial r \wedge \bar \partial r \wedge (\partial \bar \partial r)^{n-q}\}(x) =0.$ Therefore, $$\V(I^q_1(U))= \bigcup_{j=2}^N S_j,$$ which is equivalent to $$\V(I^q_1(U))-\bigcup_{j=2}^{N-1} S_j= S_N \subset M^{n+1-q}$$ since $S_i \cap S_j = \emptyset$ if $i \neq j.$ Given the definition of $M^{n+1-q},$ this means precisely that $$r_1,r_{p+2}, \dots, r_{n+1-q} \in \left(\V(I^q_1(U))-\bigcup_{j=2}^{N-1} S_j\right).$$ \qed

\section{Effective computations}
\label{effcomp}

When $x_0 \in b \Omega$ is a point of finite D'Angelo $q$-type $\Delta_q(b \Omega, x_0)=t \in \Q^+,$ we would like to give an effective upper bound for $N,$ the number of level sets of the $(n+1-q)^{th}$ commutator multitype $\ct^{n+1-q}$ in a neighborhood of $x_0,$ in terms of $t,$ the dimension $n,$ and the level of forms $q:$

\medskip
\newtheorem{effN}{Lemma}[section]
\begin{effN}
\label{effNcomp} Let $\Omega \subset \C^n$ be a smooth pseudoconvex domain, and let $x_0 \in b \Omega$ be a boundary point where $\Delta_q(b \Omega, x_0)=t.$ Let $U$ be a neighborhood of $x_0$ such that for every $x \in U \cap b \Omega$ $\Delta_q(b \Omega, x)\leq 2t^{n-q}$  and the $(n+1-q)^{th}$ commutator multitype $\ct^{n+1-q}$ takes only finitely many values $\ct^{n+1-q}_1< \dots< \ct^{n+1-q}_N$ in $U \cap b \Omega$ for some natural number $N \geq 1.$ $$N \leq (\lceil 2t^{n-q} \rceil -1)\lceil 2t^{n-q} \rceil^{\frac{(n-q)(n-q+1)}{2}-1},$$ where  $\lceil 2t^{n-q} \rceil$ is the ceiling of $2t^{n-q},$ i.e. the least integer greater than or equal to $2t^{n-q}.$
\end{effN}

\smallskip\noindent {\bf Proof:} As D'Angelo proved in \cite{opendangelo}, $\Delta_q(b \Omega, x)$ may jump in a neighborhood of $x_0.$ We can control this jump by Theorem~\ref{deltaqopenthm}, namely there exists a neighborhood $U$ of $x_0$ such that  for every $x \in U \cap b \Omega,$ $\Delta_q(b \Omega, x)\leq 2t^{n-q}.$ Shrink $U$ if necessary in order to ensure that part (i) of Theorem~\ref{cmtypethm} holds. Since the domain is pseudoconvex, by Theorem~\ref{cmtypethm} (ii) and Theorem~\ref{mlttheorem} (2), at all points of $U \cap b \Omega$ the $(n+1-q)^{th}$ commutator multitype $\ct^{n+1-q}$ takes only finitely many values $\ct^{n+1-q}_1< \dots< \ct^{n+1-q}_N$ and $\ct^{n+1-q}=(c_1, \dots, c_{n+1-q}) \in \Gamma_{n+1-q}$ with $$2 \leq c_2 \leq \dots \leq c_{n+1-q}\leq 2t^{n-q}.$$ We will now compute an upper bound for $N.$ The commutator multitype is ordered lexicographically, so even if its entries at $x_0$ are bounded above by $t$ and $\ct^{n+1-q}(x_0)=\ct^{n+1-q}_N,$ the top value, it can still happen that away from $x_0$ some entries of a lower commutator multitype are bounded above not by $t$ but by the D'Angelo type at that point, which may be higher. This is why we must work with the jump in D'Angelo type as our upper bound. In order not to have to write the quantity $2t^{n-q}$ at every step of the effective computation of $N,$ we set $t' = 2t^{n-q},$ and we will substitute back for $t'$ at the end of the argument. By Definition~\ref{wgtdef}, $c_1=1$ and the $c_i$'s are defined recursively to satisfy that $\exists \: a_1, \dots, a_k \in \N=\{0, 1, 2, \dots\}$ such that
\begin{equation}
\sum_{j=1}^k \, \frac{a_j}{c_j} =1,  \label{defsum}
\end{equation}
where the top coefficient $a_k>0.$ We know $c_1=1.$ We would like to estimate the number of possibilities for $c_2:$
$$\frac{a_1}{c_1}+ \frac{a_2}{c_2}=1.$$
If $a_1=1,$ then $\frac{a_1}{c_1}=1$ and $a_2=0,$ which contradicts the requirement that $a_2>0.$ Therefore, $a_1=0,$ and we have that $c_2=a_2.$ Now, since $2\leq c_2 \leq t' \in \Q^+,$  we estimate there are at most $\lceil t' \rceil -1$ possibilities for $a_2,$ where $\lceil t' \rceil$ is the least integer greater than or equal to $t'.$ Let us go one step further and look at $c_3.$ Now, $$\frac{a_1}{c_1}+ \frac{a_2}{c_2}+\frac{a_3}{c_3}=1.$$ As above, $a_1=0.$ There are at most $\lceil t' \rceil$ possibilities for $a_2$, namely all integers from $0$ to $\lceil t' \rceil -1$ and at most  $\lceil t' \rceil $ possibilities for $a_3,$ namely the integers from $1$ to  $\lceil t' \rceil.$ This makes a total of at most $\lceil t' \rceil^2$ possibilities for $c_3.$  If $n=3,$ we have at most  $(\lceil t' \rceil -1) \lceil t' \rceil^2$ $3$-tuples $(c_1, c_2, c_3)$ given what we got for $c_2$ and $c_3.$

Note that our method here only counts the number of possible equations for $c_j$ for $j\geq 3$ without taking into account the fact that several equations might yield the same $c_j$ or the fact that the entries $c_j$ are increasing.

In general, we expect at most $ \lceil t' \rceil^{j-1}$ possibilities for $c_j$ for $3 \leq j \leq n+1-q$ and $\lceil t' \rceil-1$ possibilities for $c_2,$ which gives at most $(\lceil t' \rceil -1) \lceil t' \rceil ^{\frac{(n-q)(n-q+1)}{2}-1}$ $(n+1-q)$-tuples $(c_1, c_2, \dots, c_{n+1-q}).$ We substitute D'Angelo type jump value for $t'$ to obtain $$N \leq (\lceil 2t^{n-q} \rceil -1)\lceil 2t^{n-q} \rceil^{\frac{(n-q)(n-q+1)}{2}-1}.$$ Clearly, this is a very generous an upper bound, which could be improved, but it will do for our purposes here. \qed

\smallskip Let us now recall from \cite{levidet} the effective bound on the vanishing order of $\text{coeff}\{\partial r \wedge \dbar r \wedge (\partial \dbar r)^{n-q}\}:$

\medskip
\newtheorem{levidetbd}[effN]{Theorem}
\begin{levidetbd}
Let $x_0 \in b \Omega$ be a point on the boundary of
the domain such that $\Delta_q(b \Omega, x_0)=t<\infty.$ At $x_0$ $\text{coeff}\{\partial r \wedge \dbar r \wedge (\partial \dbar r)^{n-q}\}$ vanishes to order at most $(\lceil t \rceil -2 )^{n-q}.$ \label{levidetbound}
\end{levidetbd}

\smallskip We would like to sharpen this result using the truncation methods from \cite{levidet} in order to derive the kind of derivative condition that appears in part (a) of the Conjectured Nullstellensatz~\ref{conjecturednullstellensatz}.

\medskip
\newtheorem{levidetderiv}[effN]{Proposition}
\begin{levidetderiv}
Let $x_0 \in b \Omega$ be a point on the boundary of
a smooth pseudoconvex domain such that $\Delta_q(b \Omega, x_0)=t<\infty$ and the rank of the Levi form at $x_0$ equals $p.$ Let $$\B_{n+1-q} (x_0) = \{r_1, r_{p+2}, \dots, r_{n+1-q}; L_2, \dots, L_{n+1-q}\}$$ be any boundary system of rank $p$ and codimension $q-1$ about the point $x_0.$ There exists a local change of variables $w=\psi(z)$ mapping $x_0$ to the origin in $\C^n$ such that the corresponding boundary system of rank $p$ about the origin in $\C^n$ in the new coordinates $$\tilde \B_{n+1-q}(0)=\{\rho_1,\rho_{p+2}, \dots, \rho_{n+1-q}; \tilde L_2, \dots, \tilde L_{n+1-q}\}$$ satisfies the following:
\begin{enumerate}[(a)]
\item $\rho_k=(\psi^{-1})^* r_k,$ the pullbacks of the functions $r_k$ under this coordinate change, are given by $\rho_k(w) = 2 Re\{w_k\} + O(|w|^2)$ in a neighborhood of the origin for $k=p+2, \dots, n+1-q;$
\item For each $k$ such that $2 \leq k \leq n+1-q,$ the coefficient of $\frac{\partial}{\partial w_k}$ in the vector field $\tilde L_k$ is non-zero at the origin;
\item For every $k$ satisfying $p+2 \leq k \leq n+1-q,$ the Levi determinant in the new coordinates $\text{coeff}\{\partial \rho \wedge \dbar \rho \wedge (\partial \dbar \rho)^{n-q}\}$ for $\rho = \rho_1 =(\psi^{-1})^*  r_1$ has a non-zero derivative at the origin of order at most $(\lceil t \rceil -2 )^{n-q}$ that involves at least one differentiation in $\frac{\partial}{\partial w_k}$ or $\frac{\partial}{\partial \bar w_k}.$
\end{enumerate}
\label{levidetderivatives}
\end{levidetderiv}

\smallskip\noindent {\bf Proof:} Translate $x_0$ to the origin in $\C^n.$ Thus, $\Delta_q(b \Omega, 0)=t<\infty.$ Let $\ct(0)=\mt(0)=(m_1, \dots, m_n)$ be the multitype at $0.$ Here we have used the pseudoconvexity of the domain $\Omega$ and part (ii) of Theorem~\ref{cmtypethm}. By part (2) of Theorem~\ref{mlttheorem}, $m_{n+1-q} \leq \Delta_q (b \Omega, 0)=t<\infty.$ If $q>1,$ it is still possible that $m_k =\infty$ for $n+2-q \leq k \leq n.$ Before we can apply Proposition~\ref{placerjprop} and Lemma~\ref{newrjlemma}, we must construct an appropriate weight $\Lambda\in \Gamma_n$ all of whose entries are finite. We distinguish two cases:

\smallskip\noindent {\bf Case 1:} $m_n < \infty.$ Then we set $\Lambda = \mt(0).$

\smallskip\noindent {\bf Case 2:} There exists $m_k =\infty$ for $n+2-q \leq k \leq n$ among the entries of $\mt(0).$ Let $k$ be the smallest integer such that $n+2-q \leq k \leq n$ and $m_k =\infty$ in $\mt(0).$ If $k=n+2-q,$ then set $\Lambda=(m_1, \dots, m_{n+1-q}, \dots, m_{n+1-q}) \in \Gamma_{n, n+1-q}$ according to Definition~\ref{elongatedef}. If $k>n+2-q,$ then set $\Lambda=(m_1, \dots, m_{n+1-q}, m_{n+2-q}, \dots, m_{k-1}, \dots, m_{k-1}) \in \Gamma_{n, k-1}.$

\smallskip\noindent Now take any boundary system $\B_{n+1-q} (0) = \{r_1, r_{p+2}, \dots, r_{n+1-q}; L_2, \dots, L_{n+1-q}\}$ of rank $p$ and codimension $q-1$ about the origin in $\C^n.$ We have that $ r = r_1 \in \mt(1; \Lambda)$ as a consequence of how we constructed $\Lambda.$ By Proposition~\ref{placerjprop}, $r_1 \in \mt(1; \Lambda)$ implies $r_k \in \mt\left(\frac{1}{\lambda_k};\Lambda\right)$ for all $k=p+2, \dots, n+1-q.$ By Lemma~\ref{newrjlemma}, there exists a coordinate change $w=\psi(z)$ in a neighborhood of the origin in $\C^n$ such that $\rho_k=(\psi^{-1})^* r_k,$ the pullbacks of the functions $r_k$ under this coordinate change for $k=p+2, \dots, n+1-q,$ satisfy $\rho_k(w) = 2 Re\{w_k\} + O(|w|^2)$ in a neighborhood of the origin. Let $\tilde \B_{n+1-q}(0)=\{\rho_1,\rho_{p+2}, \dots, \rho_{n+1-q}; \tilde L_2, \dots, \tilde L_{n+1-q}\}$ be the corresponding boundary system of rank $p$ about the origin in $\C^n$ in the new coordinates. We know that for each $k$ such that $p+2 \leq k \leq n+1-q,$ the coefficient of $\frac{\partial}{\partial w_k}$ in the vector field $\tilde L_k$ is non-zero at the origin. Furthermore, it is also evident from the special frame described on
page~\pageref{prenbd} of Section~\ref{kohnalg} and the fact that the Levi form has rank $p$ at $0$ that without loss of generality we can take $w_1, \dots, w_{p+1}$ such that $\frac{\partial \rho_1}{\partial w_1} (0) \neq 0$ and $\tilde L_2, \dots, \tilde L_{p+1}$ satisfy that the coefficient of $\frac{\partial}{\partial w_k}$ in the vector field $\tilde L_k$ is non-zero at the origin for $2 \leq k \leq n+1.$ We have already shown parts (a) and (b) of the conclusion of this proposition hold.

We just need to derive part (c). Evidently, our coordinate system $w_1, \dots, w_n$ was constructed above so that all the information about the boundary system  of rank $p$ and codimension $q-1$ about the point $x_0$ is encapsulated in variables $w_1, \dots, w_{n+1-q},$ hence we can project $\Phi: \C^n \rightarrow \C^{n+1-q}$ via $\Phi(w_1, \dots, w_n) = (w_1, \dots, w_{n+1-q}).$ Let $\tilde \rho_1$ be the push forward of the defining function $\rho=\rho_1$ under $\Phi,$ and let $\Omega'$ be the projection of $\Omega$ under $\Phi,$ i.e. the domain defined by $\tilde \rho_1.$ We claim that $\Delta_1(b \Omega', 0)=\Delta_q(b \Omega, 0)= t.$ The reason is that the $(w_1, \dots, w_{n+1-q})$ space is the image of the embedding $\phi$ that realizes the infinum in Definition~\ref{qfinitetype}. Essentially, we have constructed the special coordinates guaranteed by Theorem~\ref{mlttheorem} part (1) up to the $(n+1-q)^{th}$ coordinate, and Theorem~\ref{mlttheorem} part (2) ensures the $1$-type of the projected domain has to equal the $q$-type of the original domain. The reader should consult Catlin's original proof of the assertions in Theorem~\ref{mlttheorem} that can be found on p.555-6 of \cite{catlinbdry}. Note also that the rank of the Levi form is still $p$ at $0$ for $b\Omega'.$  As in \cite{levidet}, consider now the truncation $\widetilde {\widetilde {\rho_1}}$ of order $\lceil t \rceil$ of the Taylor expansion at $0$ of $\tilde \rho_1$ for $t =\Delta_1(b \Omega', 0)=\Delta_q(b \Omega, 0).$  As explained in \cite{levidet}, the domain $\Omega''$ defined by $\widetilde {\widetilde {\rho_1}}$ might not be pseudoconvex, but it has the same D'Angelo $1$-type $t$ at $0$ as the original domain $\Omega',$ and obviously the rank of the Levi form of $\Omega''$ at $0$ is still $p.$  Consider $\text{coeff}\{\partial \widetilde {\widetilde {\rho_1}} \wedge \dbar \widetilde {\widetilde {\rho_1}} \wedge (\partial \dbar \widetilde {\widetilde {\rho_1}})^{n-q}\},$ which is the full Levi determinant at $0$ of $b\Omega''.$ Now, let us assume there exists some $k,$ where $p+2 \leq k \leq n+1-q,$ such that $\text{coeff}\{\partial \widetilde {\widetilde {\rho_1}} \wedge \dbar \widetilde {\widetilde {\rho_1}} \wedge (\partial \dbar \widetilde {\widetilde {\rho_1}})^{n-q}\}$ is independent of both $w_k$ and $\bar w_k.$ The variety $\V\left(\left(\text{coeff}\{\partial \widetilde {\widetilde {\rho_1}} \wedge \dbar \widetilde {\widetilde {\rho_1}} \wedge (\partial \dbar \widetilde {\widetilde {\rho_1}})^{n-q}\}\right)\right)$ thus contains a complex line, which contradicts the finite type assumption on $b \Omega''.$ We have obtained the needed contradiction that shows $\text{coeff}\{\partial \widetilde {\widetilde {\rho_1}} \wedge \dbar \widetilde {\widetilde {\rho_1}} \wedge (\partial \dbar \widetilde {\widetilde {\rho_1}})^{n-q}\}$ is a polynomial of degree at most $(\lceil t \rceil -2 )^{n-q},$ which has at least one term depending on either $w_k$ or $\bar w_k$ for every $k$ such that $p+2 \leq k \leq n+1-q.$ This argument proves part (c) of Proposition~\ref{levidetderivatives} at $0 \in b \Omega''.$ We now retrace our steps. Clearly, if the conclusion of part (c) of Proposition~\ref{levidetderivatives} holds for $\text{coeff}\{\partial \widetilde {\widetilde {\rho_1}} \wedge \dbar \widetilde {\widetilde {\rho_1}} \wedge (\partial \dbar \widetilde {\widetilde {\rho_1}})^{n-q}\},$ it must also hold for the Levi determinant $\text{coeff}\{\partial \tilde \rho_1 \wedge \dbar  \tilde \rho_1 \wedge (\partial \dbar \tilde \rho_1)^{n-q}\}$ at $0 \in b \Omega'$ corresponding to the defining function $\tilde \rho_1$ before the truncation of the its Taylor expansion at $0$ took place. Furthermore, part (c) of  Proposition~\ref{levidetderivatives} must also hold for $\text{coeff}\{\partial \rho_1 \wedge \dbar  \rho_1 \wedge (\partial \dbar \rho_1)^{n-q}\}$ at $0 \in b \Omega$ since we obtained $\tilde \rho_1$ from $\rho_1$ by setting $w_{n+2-q}= \cdots = w_n =0,$ so the terms we want depending on $w_k$ or $\bar w_k$ for every $k$ in $k=p+2, \dots, n+1-q$ are present in $\text{coeff}\{\partial \rho_1 \wedge \dbar  \rho_1 \wedge (\partial \dbar \rho_1)^{n-q}\}$ as well. \qed

\section{Proof of Theorem~\ref{maintheorem}}
\label{equivpf}

We start with a natural definition that allows us to work with as small a number of level sets of the Catlin multitype as possible:

\medskip
\newtheorem{optset}{Definition}[section]
\begin{optset}
\label{optsetdef} Let $\Omega \subset \C^n$ be a smooth pseudoconvex domain, and let $x_0 \in b \Omega$ be a boundary point of finite D'Angelo $q$-type. Let $$\B_{n+1-q}(x_0)=\{r_1,r_{p+2}, \dots, r_{n+1-q}; L_2, \dots, L_{n+1-q}\}$$ be the boundary system at $x_0.$ A neighborhood $U$ of $x_0$ in $b \Omega$ is called optimal for the $(n+1-q)^{th}$ commutator multitype $\ct^{n+1-q}$ if there does not exist a smaller neighborhood $U' \subsetneq U$ such that $x_0 \in U'$ and $U'$ contains a strictly smaller number of level sets of $\ct^{n+1-q}$ than $U.$
\end{optset}

\smallskip\noindent {\bf Remark:} For any neighborhood $U$ on which the D'Angelo $q$-type is finite, the number of level sets of  the $(n+1-q)^{th}$ commutator multitype $\ct^{n+1-q}$ is finite, so $U= \bigcup_{j=1}^N S_j$ for some $N.$ We know that $x_0 \in S_N$ and $S_1$ is open  in $b \Omega.$ It follows that if $U$ is an optimal neighborhood for the $(n+1-q)^{th}$ commutator multitype $\ct^{n+1-q}$, then for every $j$ such that $2 \leq j \leq N-1,$ there exists a sequence $\{y^{(j)}_i\}_{i=1, 2, \dots}$ satisfying that  $\{y^{(j)}_i\}_{i=1, 2, \dots} \subset S_j$ and $\lim_{i \to \infty} y^{(j)}_i = x_0,$ i.e. the level sets $S_2, \dots, S_{N-1}$ accumulate at $x_0.$

\smallskip
Recall that at the first step of the modified Kohn algorithm $$\tilde I^q_1(x_0) = \sqrt[\R]{(\, r,\,
\text{coeff}\{\partial r \wedge \dbar r \wedge (\partial \dbar
r)^{n-q}\}\, )_{\ra(x_0)}}.$$ Since the defining function $r$ is identically zero on $b \Omega,$ the modified Kohn algorithm is controlled by the behavior of $\text{coeff}\{\partial r \wedge \dbar r \wedge (\partial \dbar
r)^{n-q}\}.$ This object is the Levi determinant only when $q=1;$ otherwise, $\text{coeff}\{\partial r \wedge \dbar r \wedge (\partial \dbar r)^{n-q}\}$ is a collection of ${ n-1 \choose n-q} \times { n-1 \choose n-q}$ complex-valued functions as we look at all $(n-q)^{th}$ order minors of the $(n-1) \times (n-1)$ matrix whose determinant is the Levi determinant. Let $s= { n-1 \choose n-q} \times { n-1 \choose n-q},$ and let these complex-valued Levi minors in the collection $\text{coeff}\{\partial r \wedge \dbar r \wedge (\partial \dbar r)^{n-q}\}$ be $f_1, \dots, f_s.$ For the purposes of proving Theorem~\ref{maintheorem}, we can simply consider $f=f_1 \bar f_1 + f_2 \bar f_2 + \dots +   f_s \bar f_s$ and note that $\tilde I^q_1(x_0) = \sqrt[\R]{(\, r,\, f\, )_{\ra(x_0)}}.$ The function $f$ might exhibit cancellation of derivatives, but we know it cannot be identically zero since we are assuming $\Delta_q(b \Omega, x_0)=t < \infty,$ so Theorem~\ref{levidetbound} tells us that $\text{coeff}\{\partial r \wedge \dbar r \wedge (\partial \dbar r)^{n-q}\}$ vanishes to order at most $(\lceil t \rceil -2 )^{n-q}$ at $x_0,$ hence in a neighborhood $U_{x_0}$ of $x_0$ as well, i.e. at least one of $f_1, \dots, f_s$ has a non-zero derivative on all of $U_{x_0}$ of order at most $(\lceil t \rceil -2 )^{n-q}.$ We will now prove Theorem~\ref{maintheorem}. Afterward, we will sketch how it can be strengthened subject to the Conjectured Nullstellensatz~\ref{conjecturednullstellensatz} being true and subject to obtaining slightly more information about the behavior of the sheaves of multipliers so that an effective lower bound for the subelliptic gain in the $\dbar$-Neumann problem can be computed in terms of $n,$ $t,$ and $q$ for a real-analytic pseudoconvex domain.

\medskip\noindent {\bf Proof of Theorem~\ref{maintheorem}:} The assertion that the modified Kohn algorithm finishes at step 1 densely in $b \Omega$ is a consequence of Corollary~\ref{firststepcor}. Kohn's result from \cite{kohnacta} that the modified Kohn algorithm in the real-analytic case finishes by step $2n$ is Proposition 6.20 on p.113. What we must prove here is that the number of level sets of the Catlin multitype in a neighborhood of $x_0$ acts as a counter for the Kohn algorithm.

 $\Delta_q(b \Omega, x_0)=t < \infty.$ By Theorem~\ref{deltaqopenthm}, we can shrink $U_{x_0}$ around $x_0$ to ensure that $\Delta_q(b \Omega, x) \leq 2 (\Delta_q(b \Omega, x_0))^{n-q}= 2 t^{n-q}.$ Since the D'Angelo $q$-type is finite at $x \in U_{x_0}$ and $b \Omega$ is pseudoconvex, the $(n+1-q)^{th}$ commutator multitype $\ct^{n+1-q}(x)=(c_1, \dots, c_{n+1-q})$ has only finite entries by part (2) of Theorem~\ref{mlttheorem} and part (ii) of Theorem~\ref{cmtypethm}. Therefore, a boundary system $\B_{n+1-q} (x)$ of codimension $q-1$ can be constructed at every $x \in U_{x_0}.$ It is also obvious that since the boundary is real-analytic, we can take vector fields with real-analytic coefficients in the boundary system at every $x \in U_{x_0},$ so the functions in the boundary system will also be real-analytic. Now consider the point $x_0 \in b \Omega.$ If the Levi form does not have rank at least $n-q$ at $x_0,$ let $$\B_{n+1-q}(x_0)=\{r_1,r_{p+2}, \dots, r_{n+1-q}; L_2, \dots, L_{n+1-q}\}$$ be a boundary system at $x_0$ of rank $p$ and codimension $q-1.$ We apply the Stratification Theorem, Theorem~\ref{ourstratthm}, to conclude that there exists a neighborhood $U$ of $x_0$ such that $$r_1,r_{p+2}, \dots, r_{n+1-q} \in \I \left(\V(\tilde I^q_1(U))-\bigcup_{j=2}^{N-1} S_j\right),$$ where $S_1, \dots, S_N$ are the level sets for $\ct^{n+1-q}$ in $U$ corresponding to increasing values of $\ct^{n+1-q}$ in the lexicographic ordering. We shrink $U,$ if necessary, in order to ensure the following are simultaneously satisfied:
\begin{enumerate}[(1)]
\item $U \subset U_{x_0},$ where $U_{x_0}$ is the neighborhood constructed above so that the D'Angelo type is finite and effectively bounded and $\text{coeff}\{\partial r \wedge \dbar r \wedge (\partial \dbar r)^{n-q}\}$ vanishes to order at most $(\lceil t \rceil -2 )^{n-q}$ on $U_{x_0};$
\item $U$ is optimal for the $(n+1-q)^{th}$ commutator multitype;
\item Statement~\ref{nbdprop} holds on $U$ for $k= 1, \dots, N$ in Theorem~\ref{quasiflasquesteps} (take intersection of neighborhoods guaranteed by Theorem~\ref{quasiflasquesteps} for each of the $k$ values);
\item Lemma~\ref{effNcomp} holds on $U;$
\item All parts of Proposition~\ref{stratcatlinprop} are satisfied on $U$ (in particular, the Levi form at every $x \in U$ has rank at least $p,$ its rank at $x_0.$)
\end{enumerate}
Let us look at all the points $x \in S_2 \subset U.$ We choose for each of these a neighborhood $U_x \subset U$ such that $x \in U_x$ and $U_x$ is the neighborhood guaranteed by part (i) of Proposition~\ref{stratcatlinprop}. Therefore, $U_x$ contains only two level sets $S_1$ and $S_2$ of the commutator multitype $\ct^{n+1-q}.$ Consider a boundary system $$\B_{n+1-q}(x)=\{r_1,r'_{p'+2}, \dots, r'_{n+1-q}; L'_2, \dots, L'_{n+1-q}\}$$ at $x$ defined such that on the neighborhood $U_x$ all parts of Proposition~\ref{stratcatlinprop} are satisfied (shrink $U_x,$ if necessary), where its rank $p' \geq p$ by our assumption that the rank of the Levi form is at least $p$ at every point of the big neighborhood $U.$ Note the change in notation compared to $\B_{n+1-q}(x_0).$ Since $U_x$ contains only two level sets of the Catlin multitype $\ct^{n+1-q}$ and since  and $r'_1,r'_{p'+2}, \dots, r'_{n+1-q}$ are real-analytic, $r'_1,r'_{p'+2}, \dots, r'_{n+1-q} \in \I (\V(\tilde I^q_1(U_x))).$ Recall the function $f=f_1 \bar f_1 + f_2 \bar f_2 + \dots +   f_s \bar f_s$ defined above from the Levi minors. The function $f$ is itself real-analytic, and its zero set is precisely $\V(\tilde I^q_1(U_x)).$ We can thus apply the classical \L ojasiewicz inequality to $f$ and each of the $r'_j$ functions in turn to conclude $r'_{p'+2}, \dots, r'_{n+1-q} \in \tilde I^q_1(U_x)=\sqrt[\R]{\tilde I^q_1(U_x)}.$  By Lemma~\ref{point3lemma}, $$\partial r \wedge \bar \partial r \wedge (\partial \bar \partial r)^{p'} \wedge \partial r'_{p'+2} \wedge \cdots \wedge \partial r'_{n+1-q} (x) \neq 0.$$ Consider the collection of functions $$\text{coeff}\{\partial r \wedge \bar \partial r \wedge (\partial \bar \partial r)^{p'} \wedge \partial r'_{p'+2} \wedge \cdots \wedge \partial r'_{n+1-q} \}\in \tilde I^q_2(U_x).$$ There is thus at least one function in this collection that does not vanish at the point $x.$ Let $g$ be such a function. Then $g \in \tilde I^q_2(U_x)$ and $g(x) \neq 0.$ Furthermore, there exists some neighborhood $U'_x$ such that $x \in U'_x \subset U_x$ and $g(y) \neq 0$ for every $y \in U'_x.$ Therefore, at each $x \in S_2,$ the Kohn algorithm finishes at step $2$  since we have shown there exists a non-zero subelliptic multiplier $g$ at each of those points. We conclude $\V(\tilde I^q_2(U)) \subset \bigcup_{j=3}^N S_j.$ We might no longer have strict equality as in the proof of Theorem~\ref{ourstratthm} because termination of the Kohn algorithm at a point is an open condition, so step 2 could remove not just the neighborhoods $U_x$ of each of the points $x \in S_2$ but neighborhoods of other points that may be sitting in $S_3, \dots, S_N.$ Note that $U$ was chosen so that Theorem~\ref{quasiflasquesteps} applies on it for $k=2.$ Therefore, each function $g$ that eliminates a neighborhood $U'_x$ from $\V(\tilde I^q_1(U))$ is generated by elements of $\tilde I^q_2(x_0).$ Note also that because the commutator multitype of the level set $S_2$ is strictly lower than that of $S_3, \dots, S_N,$ at least one of $r'_{p'+2}, \dots, r'_{n+1-q}$ in a boundary system at a point of $S_2$ is generated by a shorter list than at least one of the functions in the boundary system at each point of $S_3, \dots, S_N,$ which implies $$\partial r \wedge \bar \partial r \wedge (\partial \bar \partial r)^{p'} \wedge \partial r'_{p'+2} \wedge \cdots \wedge \partial r'_{n+1-q} (x) = 0 \quad \forall \, x \in S_3 \cup \cdots \cup S_N.$$ Therefore, for each $x \in S_2,$ the function $g(x)$ chosen above that eliminates the neighborhood $U'_x$ from $\V(\tilde I^q_1(U))$ vanishes on $S_3, \dots, S_N.$ Finally, we should emphasize here that our elimination of $S_2$ provides more information beyond Kohn's Proposition~\ref{stillthereprop}. The existence of a boundary system of codimension $q-1$ at every $x \in S_2$ means that $\V(\tilde I^q_1(x))$ has holomorphic dimension at most $q-1$ in the sense of Kohn at every such $x,$ so Proposition~\ref{stillthereprop} in conjunction with the \L ojasiewicz Nullstellensatz, Theorem~\ref{lojasiewicznss}, as well as Theorem~\ref{quasiflasquesteps} already guarantees that all points of $S_2$ should be eliminated at the second  step of the modified Kohn algorithm. Our method, however, explicitly constructs the elements that eliminate $S_2.$ We should note that the same will be true as we eliminate level sets $S_3, \dots, S_N$ as well.

We now have to bring the information from the level set $S_2$ forward in order to modify the zero set of $f=f_1 \bar f_1 + f_2 \bar f_2 + \dots +   f_s \bar f_s.$ For each $x \in S_2,$ we have the multiplier $g  \in \tilde I^q_2(U_x)$ satisfying that $g(y) \neq 0$ for every $y \in U'_x.$ As explained above, $g(x)=0$ on $S_3, \dots, S_N$ and $g(x)$ is generated by elements of $\tilde I^q_2(x_0).$ Now consider the collection of such functions $g(x)$ for every $x \in S_2.$ Since $\ra(x_0)$ is Noetherian, the ideal generated by all the elements of $\tilde I^q_2(x_0)$ that generate $g(x)$ for every $x \in S_2$ is finitely generated. Let $h_1, \dots, h_\beta$ be its generators. Let $f^{(2)} = h_1 \bar h_1 + \cdots + h_\beta \bar h_\beta.$ The real-analytic function $f^{(2)}\in \tilde I^q_2(x_0)$ is now real-valued and non-negative. Note that $f^{(2)} \neq 0$ for every $x \in S_2,$ $f^{(2)} \equiv 0$ on $S_3, \dots, S_N$ by construction, and $f^{(2)}$ is a subelliptic multiplier on all of $U$ since Theorem~\ref{quasiflasquesteps} holds on $U$ for $k=2.$

Next we look at the points $x \in S_3.$ We are going to recycle the notation from the previous step as all the information from $S_2$ is already being transferred via the real-analytic function $f^{(2)}.$ We choose for each of these a neighborhood $U_x \subset U$ such that $x \in U_x$ and $U_x$ is the neighborhood guaranteed by part (i) of Proposition~\ref{stratcatlinprop}. Therefore, $U_x$ contains only three level sets $S_1,$ $S_2,$ and $S_3$ of the commutator multitype $\ct^{n+1-q}.$  Consider a boundary system $$\B_{n+1-q}(x)=\{r_1,r'_{p'+2}, \dots, r'_{n+1-q}; L'_2, \dots, L'_{n+1-q}\}$$ at $x$ defined such that on the neighborhood $U_x$ all parts of Proposition~\ref{stratcatlinprop} are satisfied (up to a shrinking of $U_x$), where its rank is $p' \geq p$ just as it was the case above for $S_2.$ Since there are only three level sets of the commutator multitype in $U_x,$ it follows that Theorem~\ref{ourstratthm} applied to $U_x$ yields that $$r'_1,r'_{p+2}, \dots, r'_{n+1-q} \in \I \left(\V(\tilde I^q_1(U_x))-S_2\right),$$ where $S_1, S_2, S_3$ are the level sets for $\ct^{n+1-q}$ in $U_x.$ We may have to shrink $U_x$ a little for the previous assertion to hold. It is clear we now have to apply the \L ojasiewicz inequality to $f+ f^{(2)}$ instead of $f$ as we did at the previous step. Note that $r'_{p+2}, \dots, r'_{n+1-q}$ vanish on the zero set of $f +  f^{(2)},$ which is just the piece of $S_3$ residing in our neighborhood $U_x$ by construction. Therefore, by the \L ojasiewicz inequality, $r'_{p'+2}, \dots, r'_{n+1-q} \in \tilde I^q_2(U_x).$ By Lemma~\ref{point3lemma}, $$\partial r \wedge \bar \partial r \wedge (\partial \bar \partial r)^{p'} \wedge \partial r'_{p'+2} \wedge \cdots \wedge \partial r'_{n+1-q} (x) \neq 0,$$ where $x \in S_3.$ Consider then the collection of functions $$\text{coeff}\{\partial r \wedge \bar \partial r \wedge (\partial \bar \partial r)^{p'} \wedge \partial r'_{p'+2} \wedge \cdots \wedge \partial r'_{n+1-q} \}\in \tilde I^q_3(U_x).$$ There is thus at least one function in this collection that does not vanish at the point $x.$ Let $g$ be such a function. Then $g \in \tilde I^q_3(U_x)$ and $g(x) \neq 0.$ Furthermore, there exists some neighborhood $U'_x$ such that $x \in U'_x \subset U_x$ and $g(y) \neq 0$ for every $y \in U'_x.$ Clearly, the algorithm finishes at step 3 for every $x \in S_3.$ We now use the functions $g(x)$ for every $x \in S_3$ to construct a multiplier $f^{(3)} \in \tilde I^q_3(U)$ in the same manner we constructed $f^{(3)}$ using the Noetherian property of $\ra(x_0)$ of which $\tilde I^q_3(x_0)$ is a subideal.

Inductively, we have thus constructed $f^{(2)}, f^{(3)}, \dots, f^{(N-1)} \in \tilde I^q_{N-1}(U)$ all of which are real-analytic and eliminate $S_2, \dots, S_{N-1}.$ Let us now look back at the boundary system at $x_0$ with which we started, $$\B_{n+1-q}(x_0)=\{r_1,r_{p+2}, \dots, r_{n+1-q}; L_2, \dots, L_{n+1-q}\}.$$ The functions $r_{p+2}, \dots, r_{n+1-q}$ vanish on the zero set of $f+  f^{(2)}+  f^{(3)}+ \cdots +  f^{(N-1)},$ which is a multiplier in $\tilde I^q_{N-1}(U).$ We now apply the \L ojasiewicz inequality to $f+  f^{(2)}+  f^{(3)}+ \cdots +  f^{(N-1)}$ and each of $r_{p+2}, \dots, r_{n+1-q}$ in turn. Therefore, $r_{p+2}, \dots, r_{n+1-q} \in \tilde I^q_{N-1}(x_0).$ By Lemma~\ref{point3lemma}, $$\partial r \wedge \bar \partial r \wedge (\partial \bar \partial r)^p \wedge \partial r_{p+2} \wedge \cdots \wedge \partial r_{n+1-q} (x_0) \neq 0.$$ Furthermore, $$\partial r \wedge \bar \partial r \wedge (\partial \bar \partial r)^p \wedge \partial r_{p+2} \wedge \cdots \wedge \partial r_{n+1-q} \in \tilde  I^q_N(x_0).$$ All other points of $S_N$ in $U$ are handled in a similar manner. Clearly, the Kohn algorithm finishes by step $N$ everywhere. \qed 

%\medskip Let us move on to effectivity considerations. 

\medskip\noindent {\bf Sketch of an effective version of Theorem~\ref{maintheorem}:} Recall the real-analytic function $f=f_1 \bar f_1 + f_2 \bar f_2 + \dots +   f_s \bar f_s$ constructed above and the neighborhood $U_{x_0}$ of $x_0$ such that at least one of $f_1, \dots, f_s$ has a non-zero derivative of order at most $(\lceil t \rceil -2 )^{n-q}$ on all of $U_{x_0}.$ As mentioned above, $f$ could exhibit cancellation of derivatives, so it might not have a non-zero derivative of order at most $2(\lceil t \rceil -2 )^{n-q}$ on all of $U_{x_0}.$ As a result, we must modify its definition. Generically, we can choose real-valued polynomials $b_1, \dots, b_s$ such that $b_j(x) >0$ on $U_{x_0}$ for every $1 \leq j \leq s$ and there exists at least one non-zero derivative of order up to $4 (\lceil t \rceil -2 )^{n-q}$ (twice the bound with which we started) for  $f=b_1  f_1 \bar f_1 + b_2  f_2 \bar f_2 + \dots + b_s  f_s \bar f_s$ on all of $U_{x_0}.$ The function $f$ is real-valued and non-negative on $\C^n,$ which we can view as $\R^{2n}.$ Obviously, $\text{coeff}\{\partial r \wedge \dbar r \wedge (\partial \dbar r)^{n-q}\} = 0 \: \iff \: f=0.$ We can thus take $M = 4\, (\lceil t \rceil -2 )^{n-q}$ in the statement of the Conjectured Nullstellensatz~\ref{conjecturednullstellensatz}. 

We now need to derive condition (a) in the statement of the Conjectured Nullstellensatz~\ref{conjecturednullstellensatz} for $f$ at any $x \in U_{x_0}.$ Without loss of generality, translate $x$ to the origin. Since the domain $\Omega$ is pseudoconvex and of finite D'Angelo type, we use Proposition~\ref{levidetderivatives} to deduce that for any boundary system $$\B_{n+1-q} (0) = \{r_1, r_{p+2}, \dots, r_{n+1-q}; L_2, \dots, L_{n+1-q}\}$$ of rank $p$ at $0,$ there exists a local change of variables at $0$ such that for every $k$ satisfying $p+2 \leq k \leq n+1-q,$ $r_k(w) = 2 Re\{z_k\} + O(|z|^2)$ in a neighborhood of the origin and there is a derivative of $f$ of order at most $4 (\lceil t \rceil -2 )^{n-q}=M$ involving at least one of $L_k$ or $\bar L_k,$ which does not vanish at $0$ itself, hence in a neighborhood $V_k$ of $0.$ Note that by part (b) of Proposition~\ref{levidetderivatives} the coefficient of $\frac{\partial}{\partial z_k}$ in the vector field $L_k$ is non-zero at the origin. Since the functions $ r_{p+2}, \dots, r_{n+1-q}$ describe hypersurfaces, we can apply a real change of variables on $\C^n$ viewed as $\R^{2n}$ on a neighborhood $V_{p+2} \cap \dots \cap V_{n+1-q}$ of the origin, where the derivatives obtained from Proposition~\ref{levidetderivatives} do not vanish, so that $ r_{p+2}, \dots, r_{n+1-q}$ become 
$x_1, \dots, x_{n-q-p},$ and condition (a) stays valid. Let $U = U_{x_0} \cap V_{p+2} \dots \cap V_{n+1-q},$ $M = 4\, (\lceil t \rceil -2 )^{n-q},$ and $\mu = n-q-p.$ For each $x \in U_{x_0},$ we have constructed a neighborhood $U$ of $x$ on which condition (a) of the Conjectured Nullstellensatz~\ref{conjecturednullstellensatz} holds. Condition (b) of the Conjectured Nullstellensatz~\ref{conjecturednullstellensatz} naturally arises in the induction that proves Theorem~\ref{maintheorem} as we saw above. 

We now follow the outline of the proof of Theorem~\ref{maintheorem}. Shrink $U$ so that condition (a) of the Conjectured Nullstellensatz~\ref{conjecturednullstellensatz} holds on $U$ with respect to the boundary system at $x_0$ $$\B_{n+1-q}(x_0)=\{r_1,r_{p+2}, \dots, r_{n+1-q}; L_2, \dots, L_{n+1-q}\}.$$ This condition is added to the list of conditions on $U$ above and amounts to at most one more shrinking of the neighborhood. Now consider all the points $x \in S_2.$ Choose neighborhoods $U_x$ such that condition (a) of the Conjectured Nullstellensatz~\ref{conjecturednullstellensatz} holds. If the Conjectured Nullstellensatz~\ref{conjecturednullstellensatz} is true, then there is an effective $P$ computed from $M=4 (\lceil t \rceil -2 )^{n-q}$ and $n.$ In fact, to use the same $P$ at all levels of the induction, it is better to compute $P$ from $M=2^N (\lceil t \rceil -2 )^{n-q}$ and $n$ bearing in mind that $U$ was chosen so that Lemma~\ref{effNcomp} held on it, so we have an effective bound for $N$ in terms of  $t,$ $n,$ and $q.$ Let us use the information in Proposition~\ref{subellcostprop} to calculate the cost in loss of subelliptic gain in eliminating $S_2.$ Parts (ii) and (vi) of Proposition~\ref{subellcostprop} imply $f$ is a multiplier with $\epsilon= \frac{1}{2}.$ Capturing $r'_{p'+2}, \dots, r'_{n+1-q}$ via the Conjectured Nullstellensatz~\ref{conjecturednullstellensatz} comes at a cost of dividing the gain by at most $P,$ so we have $\epsilon\geq \frac{1}{2P}$ for each of them by part (iv) of Proposition~\ref{subellcostprop}. Note that the $M$ is effectively computed in terms of $t,$ $n,$ and $q$ and universal on $U,$ so $P$ will also be a function of $M$ and $n,$ hence of $t,$ $n,$ and $q.$ The application of Lemma~\ref{point3lemma}, costs another factor of $\frac{1}{2},$ so complex gradients $\partial r'_{p'+2}, \dots, \partial r'_{n+1-q}$ correspond to a gain of $\epsilon\geq \frac{1}{4P}$ by part (v) of Proposition~\ref{subellcostprop}. Finally, taking the determinant $\partial r \wedge \bar \partial r \wedge (\partial \bar \partial r)^{p'} \wedge \partial r'_{p'+2} \wedge \cdots \wedge \partial r'_{n+1-q}$ still leaves us with a gain of $\epsilon\geq\frac{1}{4P}$ by part (vi) of Proposition~\ref{subellcostprop}. Therefore, for every point $x \in S_2$ and the corresponding $g(x),$ we have the lower bound for subelliptic gain $\epsilon\geq \frac{1}{4P}$ that works on the neighborhood $U_x.$ Now we hit the next significant issue that was listed in the introduction. While we know from Theorem~\ref{quasiflasquesteps} that $g(x)$ is generated by elements of $\tilde I^q_2(x_0),$ Theorem~\ref{quasiflasquesteps} is a qualitative result. It does not tell us to what subelliptic gain those elements that generate $g(x)$ correspond in a neighborhood of $x_0$ or even better on $U.$ In other words, we would need a quantitative version of Theorem~\ref{quasiflasquesteps} in order to compute the subelliptic gain corresponding to $f^{(2)}$ even if it involved a shrinking of the neighborhood $U.$ After all, $U$ was chosen to be optimal, so no shrinking of it can diminish the number $N$ of level sets of the Catlin multitype. 

If such a quantitative version of Theorem~\ref{quasiflasquesteps} could be proven, then there would be no further roadblocks to an effective computation of subelliptic gain. Assume that $f^{(2)}$ corresponded to subelliptic gain $P_2$ effectively computed, then we would consider all points $x \in S_3$ and take generic real-valued polynomials $\tilde b_1$ and $\tilde b_2$ such that $\tilde b_j(x) >0$ on $U$ and there exists at least one non-zero derivative of order up to $8 (\lceil t \rceil -2 )^{n-q}$ for $\tilde b_1 f + \tilde b_2 f^{(2)}$ on $U.$ Neighborhoods $U_x$ would thus be chosen so that condition (a) in the Conjectured Nullstellensatz~\ref{conjecturednullstellensatz} holds. The sum $\tilde b_1 f + \tilde b_2 f^{(2)}$ would be a subelliptic multiplier with gain bounded below by the minimum of the gain for each of the pieces, i.e. $\epsilon\geq \min\{ \frac{1}{4P}, P_2\}.$ As a result of the Nullstellensatz, $r'_{p'+2}, \dots, r'_{n+1-q}$ would then correspond to subelliptic gain satisfying $\epsilon\geq \min\{ \frac{1}{4P^2},\frac{P_2}{P}\}.$ The subelliptic gain of the functions $g(x)$ that eliminate $S_3$ would then be $\epsilon\geq \min \{ \frac{1}{8P^2},\frac{P_2}{2P}\}.$ Continuing the process under the assumption a quantitative version of Theorem~\ref{quasiflasquesteps} could be proven, if $f^{(3)}, \dots, f^{(N-1)}$ came with subelliptic gain $P_3, \dots, P_{N-1}$ respectively, then in the end we would get $\epsilon \geq \min\{ \frac{1}{2 (2P)^{N-1}},\frac{P_2}{(2P)^{N-2}}, \dots, \frac{P_{N-1}}{2P} \}.$ \qed

\bibliographystyle{plain}
\bibliography{SmoothTypeEquiv}

\begin{thebibliography}{10}

\bibitem{bochnakrisler}
J.~Bochnak and J.~J. Risler.
\newblock Sur les exposants de {L}ojasiewicz.
\newblock {\em Comment. Math. Helv.}, 50(4):493--507, 1975.

\bibitem{bazilandreea}
Vasile Brinzanescu and Andreea~C. Nicoara.
\newblock On the relationship between d'{A}ngelo q-type and {C}atlin q-type.
\newblock \textsl{{P}reprint} ar{X}iv:1302.2294v4, [math.CV] 8 Jan 2014,
  {D}{O}{I}: 10.1007/s12220-014-9490-5; in press at the {J}ournal of
  {G}eometric {A}nalysis.

\bibitem{catlinnec}
David Catlin.
\newblock Necessary conditions for subellipticity of the {$\bar \partial
  $}-{N}eumann problem.
\newblock {\em Ann. of Math. (2)}, 117(1):147--171, 1983.

\bibitem{catlinbdry}
David Catlin.
\newblock Boundary invariants of pseudoconvex domains.
\newblock {\em Ann. of Math. (2)}, 120(3):529--586, 1984.

\bibitem{catlinsubell}
David Catlin.
\newblock Subelliptic estimates for the {$\overline\partial$}-{N}eumann problem
  on pseudoconvex domains.
\newblock {\em Ann. of Math. (2)}, 126(1):131--191, 1987.

\bibitem{catlincho}
David~W. Catlin and Jae-Seong Cho.
\newblock Sharp estimates for the $\overline\partial$-neumann problem on
  regular coordinate domains.
\newblock \textsl{{P}reprint.} arXiv:0811.0830v1, [math.CV] 5 Nov 2008.

\bibitem{fribourgcatda}
David~W. Catlin and John~P. D'Angelo.
\newblock Subelliptic estimates.
\newblock In {\em Complex analysis}, Trends Math., pages 75--94.
  Birkh\"auser/Springer Basel AG, Basel, 2010.

\bibitem{opendangelo}
John~P. D'Angelo.
\newblock Real hypersurfaces, orders of contact, and applications.
\newblock {\em Ann. of Math. (2)}, 115(3):615--637, 1982.

\bibitem{dangelo}
John~P. D'Angelo.
\newblock {\em Several complex variables and the geometry of real
  hypersurfaces}.
\newblock Studies in Advanced Mathematics. CRC Press, Boca Raton, FL, 1993.

\bibitem{df}
Klas Diederich and John~E. Fornaess.
\newblock Pseudoconvex domains with real-analytic boundary.
\newblock {\em Ann. Math. (2)}, 107(2):371--384, 1978.

\bibitem{khanhzampieri}
Tran~Vu Khanh and Giuseppe Zampieri.
\newblock Precise subelliptic estimates for a class of special domains.
\newblock \textsl{{P}reprint.} arXiv:0812.2560v2, [math.CV] 7 Jan 2009.

\bibitem{dbneumann1}
J.~J. Kohn.
\newblock Harmonic integrals on strongly pseudo-convex manifolds. {I}.
\newblock {\em Ann. of Math. (2)}, 78:112--148, 1963.

\bibitem{dbneumann2}
J.~J. Kohn.
\newblock Harmonic integrals on strongly pseudo-convex manifolds. {II}.
\newblock {\em Ann. of Math. (2)}, 79:450--472, 1964.

\bibitem{kohnacta}
J.~J. Kohn.
\newblock Subellipticity of the {$\bar \partial $}-{N}eumann problem on
  pseudo-convex domains: sufficient conditions.
\newblock {\em Acta Math.}, 142(1-2):79--122, 1979.

\bibitem{kollareffpol}
J{\'a}nos Koll{\'a}r.
\newblock An effective \l ojasiewicz inequality for real polynomials.
\newblock {\em Period. Math. Hungar.}, 38(3):213--221, 1999.

\bibitem{classicalloj}
S.~{\L}ojasiewicz.
\newblock Sur le probl\`eme de la division.
\newblock {\em Studia Math.}, 18:87--136, 1959.

\bibitem{lojasiewicz}
Stanislas {\L}ojasiewicz.
\newblock Sur la g\'eom\'etrie semi- et sous-analytique.
\newblock {\em Ann. Inst. Fourier (Grenoble)}, 43(5):1575--1595, 1993.

\bibitem{andreeaqf}
Andreea~C. Nicoara.
\newblock Coherence and other properties of sheaves in the {K}ohn algorithm.
\newblock \textsl{{P}reprint.} arXiv:1308.5289v1, [math.AG] 24 Aug 2013,
  {D}{O}{I}: 10.1142/{S}0129167{X}14500773; in press at the {I}nternational
  {J}ournal of {M}athematics.

\bibitem{levidet}
Andreea~C. Nicoara.
\newblock Effective vanishing order of the {L}evi determinant.
\newblock {\em Math. Ann.}, 354(4):1223--1245, 2012.

\bibitem{siunote}
Yum-Tong Siu.
\newblock Effective termination of {K}ohn's algorithm for subelliptic
  multipliers.
\newblock {\em Pure Appl. Math. Q.}, 6(4, Special Issue: In honor of Joseph J.
  Kohn. Part 2):1169--1241, 2010.

\bibitem{solerno}
Pablo Solern{\'o}.
\newblock Effective \l ojasiewicz inequalities in semialgebraic geometry.
\newblock {\em Appl. Algebra Engrg. Comm. Comput.}, 2(1):2--14, 1991.

\bibitem{thom}
Ren{\'e} Thom.
\newblock On some ideals of differentiable functions.
\newblock {\em J. Math. Soc. Japan}, 19:255--259, 1967.

\bibitem{tougeronqf}
Jean-Claude Tougeron.
\newblock Faisceaux diff\'erentiables quasi-flasques.
\newblock {\em C. R. Acad. Sci. Paris}, 260:2971--2973, 1965.

\bibitem{tougeron}
Jean-Claude Tougeron.
\newblock {\em Id\'eaux de fonctions diff\'erentiables}.
\newblock Springer-Verlag, Berlin, 1972.
\newblock Ergebnisse der Mathematik und ihrer Grenzgebiete, Band 71.

\end{thebibliography}

\end{document}